\documentclass{article}%
\usepackage{eurosym}
\usepackage{amsmath}
\usepackage{amssymb}
\usepackage{amsfonts}
\usepackage{geometry}
\usepackage{graphicx}
\usepackage{float}
\usepackage{caption}
\usepackage{cite}
\usepackage{subcaption}
\usepackage{verbatim}
\usepackage[utf8]{inputenc}

\providecommand{\U}[1]{\protect\rule{.1in}{.1in}}
\newtheorem{theorem}{Theorem}

\newtheorem{definition}[theorem]{Definition}
\newtheorem{example}[theorem]{Example}

\newtheorem{lemma}[theorem]{Lemma}
\newtheorem{notation}[theorem]{Notation}

\newtheorem{remark}[theorem]{Remark}

\newenvironment{proof}[1][Proof]{\noindent\textbf{#1.} }{\ \rule{0.5em}{0.5em}}
\begin{document}

\title{On the Centred Hausdorff Measure of the Sierpinski Gasket}

\author{Marta LLorente$^1$, Mar\'{i}a Eugenia Mera$^2$ and Manuel Mor\'{a}n $^{2,3}$}
\date { }
\maketitle

{\centering{\small {$^{1}$ Departamento de An\'{a}lisis Econ\'{o}mico: Econom%
\'{\i}a Cuantitativa. Universidad Aut\'{o}noma de Madrid, Campus de
Cantoblanco, 28049 Madrid, Spain.}}}
\newline
{\centering{\small {$^{2}$ Departamento de An\'{a}lisis Econ\'{o}mico y
Econom\'{\i}a Cuantitativa. Universidad Complutense de Madrid. Campus de
Somosaguas, 28223 Madrid, Spain.}}}
\newline
{\centering{{\small {$^{3}$ IMI-Institute of Interdisciplinary Mathematics.
Universidad Complutense de Madrid. Plaza de Ciencias 3, 28040 Madrid, Spain.}%
}}}

{\centering{Emails: m.llorente@uam.es, mera@ucm.es, mmoranca@ucm.es}}

\textit{Short Title:} On the Centred Hausdorff Measure of the Sierpinski Gasket

\begin{abstract}
We show that the centred Hausdorff measure, $C^{s}(S),$ with $s=\frac{\log3}{\log2},$ of the Sierpinski gasket $S$, is $C$-computable
(continuous-computable), in the sense that its value is the solution of the minimisation problem of a continuous function on a compact domain. We also show that $C^{s}(S)$ is $A$-computable (algorithmic-computable) in the sense that there is an algorithm that converges to $C^{s}(S),$ with error bounds tending to zero. Using this algorithm and bounds we show that $C^{s}(S)\sim1.0049,$ and we establish a conjecture for the value of the spherical Hausdorff $s$-measure of $S,$ $\mathcal{H}_{sph}^{s}(S)\sim0.8616$, and
provide an upper bound for it, $\mathcal{H}_{sph}^{s}(S)\leq0.8619.$
\end{abstract}

\textit{Keywords}: Self-Similar Sets, Sierpinski Gasket, Hausdorff Measures, Density of Measures, Computability of Fractal Measures.
\newline
\textit{[2020] MSC: 28A78, 28A80, 28A75}
\newpage


\section{Introduction and main results\label{Section Introduction}}

The role played by the Sierpinski gasket in the start of Fractal Geometry
might well be compared to that played by the triangle or the circle in the
start of Euclidean Geometry. The vast amount of literature devoted to this
iconic fractal gives account of its relevance (for instance, more than one
hundred articles only in mathscinet database in the last three years). In
spite of this fact some basic parameters of the Sierpinski gasket remain
unknown. The complex nature of fractal objects has given rise to a great
variety of measurement tools such as metric measures and dimensions in
Euclidean spaces. The Hausdorff measure, $\mathcal{H}^{s},$ the centred
Hausdorff measure, $C^{s},$ the spherical Hausdorff measure, $\mathcal{H}%
_{sph}^{s},$ and the packing measure, $P^{s},$ (see Sec.~\ref{subsection
metric measures} for definitions) are some of the existing ones that are so
necessary to capture the multiple geometric properties of fractal objects.
In these measures, which we shall refer to as $s$-measures, the exponent $s$
is a real positive number. We shall denote by $\mathcal{M}^{s}$ the set $%
\left\{ \mathcal{H}^{s},C^{s},P^{s},\mathcal{H}_{sph}^{s}\right\} $ of
metric $s$-measures.

For integer values of $s$ all metric $s$-measures are multiples of the
corresponding Lebesgue measure $\mathcal{L}^{s},$ so that they are well
understood and computable for some families of conspicuous subsets.

If $s$ is not an integer, the examples of $s$-sets (i.e. sets with finite
and positive $s$-measure), for which the exact value of any of these
measures is known, may be considered exceptional even in the class of
self-similar sets with strong separation condition properties (see, for
example, \cite{AYER, BZ, CMP, FENG2, GZ, Me, Z1} and the references
therein). In fact, the exact value of the mentioned $s$-measures of almost
any self-similar set, including the Sierpinski triangle or the Koch curve,
is not known.

The nature of the results shown in this paper for the centred Hausdorff
measure (see also \cite{LLMM2} for an example with the packing measure)
indicates that finding the exact value of a metric measure of a self-similar
set is essentially a problem of a computational nature. As a rule of thumb,
the required computations are easier when the constituent parts of the
corresponding set are suitably separated, but as this separation becomes
shorter, the study of their properties becomes computationally more arduous.

In this paper we undertake the issue of computing the centred Hausdorff
measure of the Sierpinski gasket, $S,$ a self-similar set for which the
separation among its constituent parts is zero, but it satisfies the open
set condition (OSC for the sequel, see Sec. \ref{siersect} below). To the
best of our knowledge, this is the first known computation of the centred
Hausdorff measure of a connected self-similar set with OSC. In Theorem \ref%
{Reduction1} we show that $C^{s}(S),$ with $s=\frac{\log 3}{\log 2},$ is $C$%
-computable, in the sense that its value is the solution of a minimisation
problem on a compact domain. See \cite{AYER} for an example of a
self-similar set $K$ in the line, with OSC, for which the Hausdorff measure
is not $C$-computable. In Theorem \ref{theorem bounds discrete} we show that 
$C^{s}(S)$ is $A$-computable, providing an algorithm whose output converges
to $C^{s}(S),$ with error bounds tending to zero.

The many efforts directed towards the problem of the computability of metric
measures of self-similar sets have yielded interesting results on the nature
of the optimal coverings necessary for the calculation of the Hausdorff
measure of $S,$ and even quite precise estimates of this value (see \cite%
{MARI, Mora, J, JZZ, ZF}, and the references therein). Using the method
introduced in \cite{Mo1} we shall establish a conjecture for the spherical
Hausdorff measure of $S$ and discuss the relationship with the results
mentioned above. We shall first summarise some basic definitions and
notation to understand the problem.

\subsection{The Sierpinski gasket\label{siersect}}

\begin{figure}[H]
\centering
\includegraphics[scale=0.5]{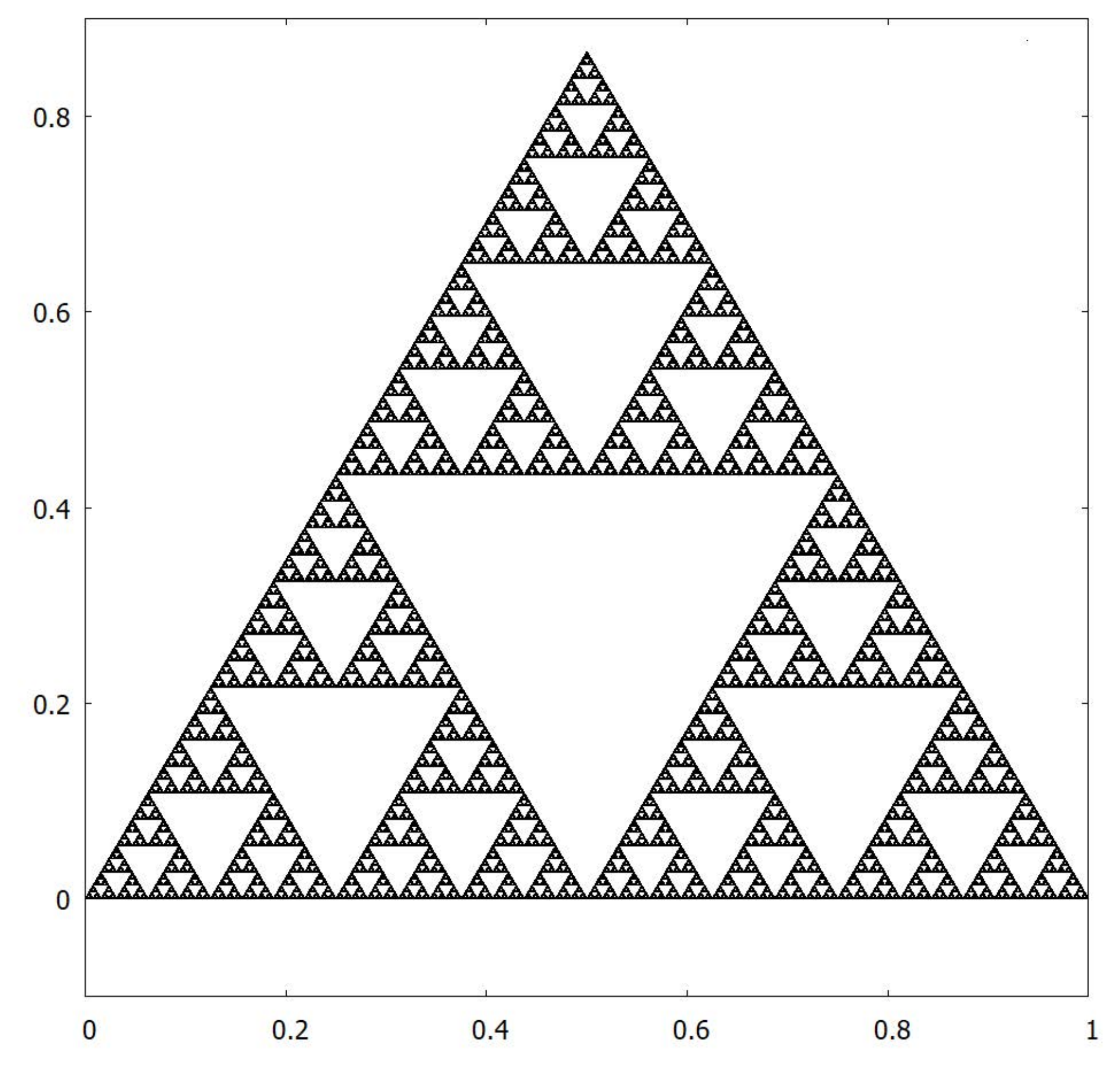}
\caption{Sierpinski gasket $S.$ }
\label{sierpinski}
\end{figure}

The Sierpinski gasket or Sierpinski triangle (see Fig. ~\ref{sierpinski}) is
a special case of a self-similar set which is generated by a system $\Psi
=\{f_{0,}f_{1,}f_{2}\}$ of three contracting similitudes of the plane, with
contraction ratios $r_{i}:=\frac{1}{2},$ $i\in M:=\{0,1,2\},$ given by 
\begin{equation*}
f_{i}(x)=\frac{1}{2}x+v_{i},\text{ }i\in M\text{ }\ \text{where }v_{0}=(0,0),
\text{ }v_{1}=(\frac{1}{2},0),\text{ and }v_{2}=\frac{1}{4}(1,\sqrt{3}).
\end{equation*}

We use composite indices, $i:=i_{1},i_{2},...,i_{k}\in M^{k},$ to denote the
compositions $f_{i}:=f_{i_{1}}\circ f_{i_{2}}\circ ...\circ f_{i_{k}}$ and
we write $r_{i}$ for the contraction ratio of $f_{i}$ (which equals $2^{-k}$
if $i\in M^{k}).$ We shall denote as $z_{i}$ to the fixed point of each $%
f_{i},$ $i\in M,$ that is, $z_{i}=2v_{i},$ $i\in M.$

The Sierpinski gasket $S,$ as the attractor of $\Psi,$ is the invariant set
of the Hutchinson operator, $F,$ defined, for $A\subset\mathbb{R}^{2},$ by 
\begin{equation}
F(A):=f_{0}(A)\cup f_{1}(A)\cup f_{2}(A),  \label{hutchin}
\end{equation}
$S$ being the unique non-empty compact set admitting the self-similar
decomposition 
\begin{equation*}
S=f_{0}(S)\cup f_{1}(S)\cup f_{2}(S)=F(S).
\end{equation*}

$S$ can be parameterised as $S=\left\{ \pi(i):i\in\Sigma\right\} $ with
parameter space $\ \ \Sigma:=M^{\infty}$ and \textit{geometric projection
mapping} $\pi:\Sigma\rightarrow S$ given by $\pi(i)=\cap_{k=1}^{%
\infty}f_{i(k)}(S),\ $ where $i(k)\in M^{k}$ denotes the $k$-th curtailment $%
i_{1}\dots i_{k}$ of $i=i_{1}i_{2}\dots\in\Sigma.$ Notice that $\pi$ is
non-injective. We adopt the convention $M^{0}=\emptyset$ and write $%
M^{\ast}=\cup_{k=0}^{\infty}M^{k}$ for the set of words of finite length.
For any $i\in M^{\ast},$ the cylinder sets are denoted by $S_{i}:=f_{i}(S),$
and $S_{i}:=S,$ if $i\in M^{0}.$ For $i\in M^{k},$ $S_{i}$ is a \textit{%
cylinder of the }$k$-th generation, or $k$-cylinder.

\begin{figure}[H]
\centering
\includegraphics[scale=0.5]{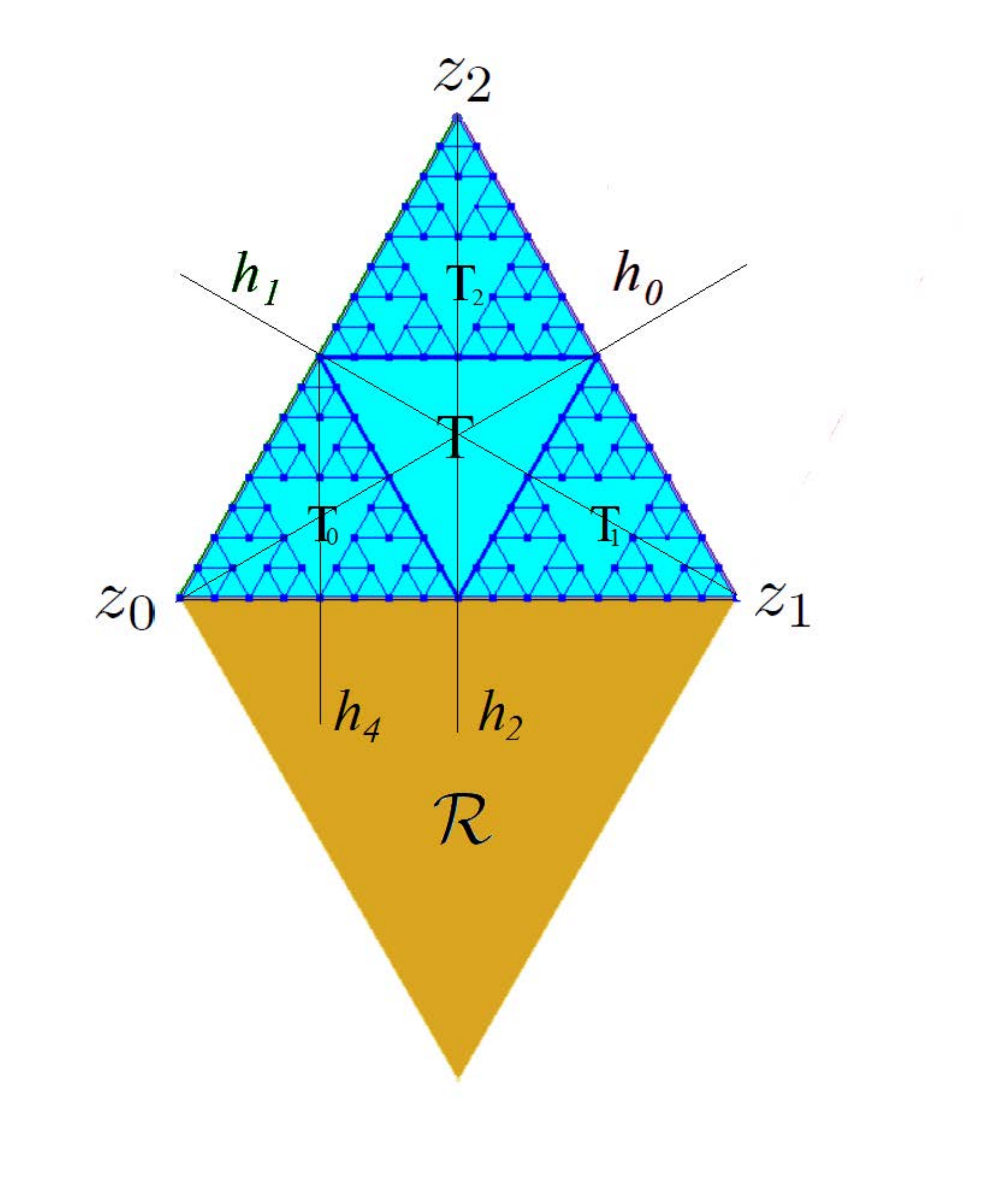}
\caption{A feasible open set $\mathcal{R}$ for $S.$\newline
The rhombus $\mathcal{R}$ is the topological interior of the union of $T$
and its reflection across the edge opposite the point $z_{2}.$ Triangles $T$
and $T_{i},$ $i\in M.$ Altitudes $h_{i}$ of $T$ through $z_{i},$ $i\in M,$
and altitude $h_{4}$ of $T_{0}$ through $f_{2}(z_{0}).$}
\label{Feasibleopenset}
\end{figure}

\begin{notation}
\label{convexhull}Throughout the document the equilateral triangle whose
vertices are $z_{i},$ $i\in M$ will be denoted by $T,$ and we shall
abbreviate $T_{i}:=f_{i}(T)$ if $i$ $\in M^{\ast },$ $k\in \mathbb{N}^{+},$
and $T_{i}:=T$ if $i\in M^{0}$ (see Fig.~\ref{Feasibleopenset}).
\end{notation}

The system of similitudes $\Psi $ satisfies the OSC (see \cite{MoP}) meaning
that there is an open set $\mathcal{R\subset }\mathbb{R}^{2}$ satisfying $%
f_{i}\mathcal{(R)\subset R}$ for all $i\in M$ and $f_{i}(\mathcal{R)\cap }%
f_{j}(\mathcal{R)=\varnothing }$ for $i,$ $j\in M,$ $i\neq j.$ We will refer
to such a set $\mathcal{R}$ as \textit{a} \textit{feasible open set} (for $%
S).$ Furthermore, if $\mathcal{R\cap }S\neq \varnothing $ then $\mathcal{R}$
satisfies the strong open set condition SOSC (cf. \cite{La, Sc, Mo2}). One
feasible open set that fulfils the SOSC is the open rhombus $\mathcal{R}$
composed of the topological interior of the union of $T$ and its reflection
across the edge of $T$ opposite the point $z_{2}$ (see Fig.~\ref{Feasibleopenset}). As $\Psi $ satisfies the OSC, the dimension of all
metric measures are the same (see \cite{Falc, Mat1} for the notion of
dimension of a measure), and they also coincide with the similarity
dimension,$\ \dim S=\frac{\log 3}{\log 2},$ which is the value that
satisfies $\sum_{i=0}^{2}r_{i}^{\dim S}=1.$

The Sierpinski gasket $S$ can also be considered of as a probability measure
supported by $S.$ Let $\mathcal{P}(\mathbb{R}^{2})$ be the space of
compactly supported probability Borel measures on $\mathbb{R}^{2}$ and let $%
\mathbf{M:}$ $\mathcal{P}(\mathbb{R}^{2})\rightarrow \mathcal{P}(\mathbb{R}%
^{2})$ be the Markov operator, defined by%
\begin{equation*}
\mathbf{M}(\alpha )=\sum_{i=0}^{2}\frac{1}{3}\left( \alpha \circ
f_{i}^{-1}\right) ,\text{ }\alpha \in \mathcal{P}(\mathbb{R}^{2}).
\end{equation*}%
The operator $\mathbf{M}$ is contractive on $\mathcal{P}(\mathbb{R}^{2}),$
equipped with a suitable metric (see \cite{HUTCH, BARNS}). Its unique fixed
point, $\mu ^{\ast },$ is called the \textit{invariant} or, sometimes, 
\textit{natural probability measure.} It is supported on $S$ and satisfies 
\begin{equation}
\mathbf{M}^{k}(\alpha )=3^{-k}\sum_{i\in M^{k}}\left( \alpha \circ
f_{i}^{-1}\right) \xrightarrow[k \to \infty]{w}\mu ^{\ast }
\label{invariant}
\end{equation}%
for any $\alpha \in \mathcal{P}(\mathbb{R}^{2}).$ Here, $\overset{w}{%
\rightarrow }$ denotes the weak convergence of measures and $\mathbf{M}^{k}=%
\mathbf{M}\circ \mathbf{M}\circ ...\circ \mathbf{M}$ is the $k$-th iterate
of $\mathbf{M.}$ Furthermore, $\mu ^{\ast }$ coincides with the projection
on $S$ of $\nu ,$ 
\begin{equation}
\mu ^{\ast }=\nu \circ \pi ^{-1},  \label{projection}
\end{equation}%
where $\nu $ is the Bernoulli measure on $\Sigma $ that gives weight $\frac{1%
}{3}$ to each symbol in $M$ (see \cite{HUTCH}).

By \eqref{projection}, we know that $\mu^{\ast}(S_{i})=3^{-k}$ for $i\in
M^{k}.$ On the other hand, any metric measure $\alpha\in\mathcal{M}^{\dim S}$
scales under similitudes, i.e. $\alpha(S_{i})=r_{i}^{\dim
S}\alpha(S)=3^{-k}\alpha(S)$, for $i\in M^{k}.$ Since $\mu^{\ast}$ and $%
\alpha,$ are multiples on cylinder sets, they are indeed multiple measures,
and then, all the normalised measures $\left( \alpha(S)\right)
^{-1}\alpha\lfloor_{S}$ coincide with $\mu^{\ast}$ and with the normalised
Hausdorff measure $\mu:=\frac{\mathcal{H}^{s}\lfloor_{S}}{\mathcal{H}^{s}(S)}%
,$ with $s=\dim S.$ Here $\beta\lfloor_{S}$ stands for the restriction of
the measure $\beta$ to $S.$

\begin{remark}
\label{invariant measure} From now on we shall work with the measure $\mu $
rather than with other metric measures on $S,$ bearing in mind that if $%
s=\dim S,$ $A$ is a Borel set and $\alpha \in \mathcal{M}^{s}\mathcal{%
\lfloor }S\mathcal{\,}$\ with $\mathcal{M}^{s}\mathcal{\lfloor }S:=\left\{ 
\mathcal{H}^{s}\lfloor _{S},\mathcal{H}_{sph}^{s}\lfloor _{S},C^{s}\lfloor
_{S},P^{s}\lfloor _{S}\right\} ,$ we have 
\begin{equation}
\alpha (A)=\alpha (S)\mu (A).  \label{muvshauss}
\end{equation}%
Whence the computation of $\alpha (A)$ boils down to the computation of $%
\alpha (S)$ plus the computation of $\mu (A)$ (see an example of the
estimation of $C^{s}\lfloor _{S}(A)$ when $A$ is a ball in Example \ref%
{cs-of-ball} at the end of Sec. \ref{section spherical}).
\end{remark}

We now introduce the following notation.

\begin{notation}
\textbf{\ }\label{typicalball} We write $\mathcal{B}$ for the set of closed
balls $B(x,d)$ centred at $x\in S$ and with radius $d>0$ satisfying that
there is some feasible open set $\mathcal{R}$ for $S$ with $B(x,d)\subset%
\mathcal{R}. $
\end{notation}


\subsection{Metric measures and their relationship with s-densities\label%
{subsection metric measures}}

Following Saint Raymond and Tricot (see \cite{SRT}), the Hausdorff centred
measure, $C^{s}(A),$ of a subset $A\subset\mathbb{R}^{n}$ is defined in a
two-step process (see \cite{EDG} for definitions on general metric spaces).
First, the premeasure $C_{0}^{s}(A)$ is defined for any $s>0$ by 
\begin{equation}
C_{0}^{s}(A)=\lim_{\delta\rightarrow0}\inf\left\{ \sum\limits_{i=1}^{\infty
}(2d_{i})^{s}\ :\ \ 2d_{i}\leq\delta,\text{ }i=1,2,\dots\right\} ,
\label{premeasure}
\end{equation}
where the infimum is taken over all coverings, $\left\{
B(x_{i},d_{i})\right\}_{i\in\mathbb{N}^{+}},$ of $A$ by closed balls $%
B(x_{i},d_{i})$ centred at points $x_{i}\in A.$ The second step, needed by
the lack of monotonicity of $C_{0}^{s}(A)$ due to the restriction $x_{i}\in
A,$ $i\in\mathbb{N}^{+}$ (cf. \cite{TO2} and \cite[Example~4 ]{LLM2}), leads
to the following definition of the \textit{centred Hausdorff }$s$\textit{%
-dimensional measure } 
\begin{equation*}
C^{s}(A)=\sup\left\{ C_{0}^{s}(F):\text{ \ }F\subset A,\ F\text{ closed}%
\right\}.
\end{equation*}

However, in \cite{LLM2} it is proved that, if $E$ is a self-similar set with
OSC, $A\subset\mathbb{R}^{n}$ is a compact set, and $s$ is its similarity
dimension \cite{HUTCH}$,$ then $C^{s}(A\cap E)=C_{0}^{s}\left( A\cap
E\right) $ implying that the second step can be omitted. As we shall see,
this makes it possible to reduce the problem of calculating this fractal
measure to the computation of the optima of certain density functions (see %
\eqref{density}).

In regard to metric measures based on packings, the standard packing
measure, $P^{s},$ is relevant in this research. It is defined in a two steps
process, 
\begin{equation*}
P_{0}^{s}(A)=\lim_{\delta\rightarrow0}\sup\left\{ \sum\limits_{i=1}^{\infty
}(2d_{i})^{s}:\text{ \ \ }2d_{i}\leq\delta,\text{ }i=1,2,\dots\right\},
\end{equation*}
where the supremum above is taken over all $\delta$-\textit{packings} $%
\left\{ B(x_{i},d_{i})\right\} _{i\in\mathbb{N}^{+}},$ with $x_{i}\in A$ for
all $i,$ and $B(x_{i},d_{i})\cap B(x_{j},d_{j})=\varnothing$ for $i\neq j,$
and 
\begin{equation*}
P^{s}(A)=\inf\left\{ \sum\limits_{i=1}^{\infty}P_{0}^{s}(F_{i})\right\},
\end{equation*}
where the infimum is taken over all coverings $\left\{ F_{i}\right\} _{i\in%
\mathbb{N}^{+}}$ of $A$ by closed sets $F_{i}.$ When $A$ is a compact set
with $P_{0}^{s}(A)<\infty,$ then $P^{s}(A)=P_{0}^{s}(A)$ holds, and the
second step above might be omitted (see \cite{FENG}).

Sec.~\ref{section spherical} is devoted to the \textit{spherical} $s$\textit{%
-dimensional Hausdorff measure, }$\mathcal{H}_{sph}^{s}(A),$ which is
obtained by removing in \eqref{premeasure} the requirement that the balls
are centred at points of $A.$ Finally, the classical $s$\textit{-dimensional
Hausdorff measure, }$\mathcal{H}^{s}(A),$ results if coverings of $A$ by
arbitrary subsets, $\left\{ U_{i}\right\} _{i\in \mathbb{N}^{+}},$ are
considered and $2d_{i}$ is replaced in \eqref{premeasure} with the diameter
of $U_{i},$ $\left\vert U_{i}\right\vert .$ No second step is required for
these two last measures.

Let us recall that, for $\alpha\in\cup_{s>0}\mathcal{M}^{s}$ $\ $and $%
A\subset\mathbb{R}^{n},$ the $\alpha$-dimension of $A$ is defined by $%
\dim_{\alpha }A=\inf\left\{ s:\alpha(A)=0\right\} $ and that $%
\dim_{\alpha}(S)=\frac{\log3}{\log2}$ (see \cite{Mat1} for the definition of
the dimension of a set with respect to a metric measure).

We now introduce the $s$-densities, a useful tool for analysing the
behaviour of metric measures defined for measures $\alpha\in\mathcal{M}%
^{s}\lfloor_{A}$ by 
\begin{equation}
\theta_{\alpha}^{s}(x,d)=\frac{\alpha(B(x,d))}{(2d)^{s}},\ s>0,\text{ }x\in%
\mathbb{R}^{n},\ d>0.  \label{density}
\end{equation}

A fundamental result in geometric measure theory (Marstrand's theorem, \cite%
{MARST}) states that, except when $s$ is an integer, the limiting values 
\begin{equation*}
\overline{\theta }_{\alpha }^{s}(x)=\limsup_{d\rightarrow 0}\theta _{\alpha
}^{s}(x,d)\text{ and \ }\underline{\theta }_{\alpha
}^{s}(x)=\liminf_{d\rightarrow 0}\theta _{\alpha }^{s}(x,d),
\end{equation*}%
called \textit{upper and lower }$s$-densities respectively, cannot coincide
on subsets $A\subset \mathbb{R}^{n}$ of positive $\alpha $-measure.
Furthermore, classical results in fractal geometry (see, for instance, \cite%
{MAT2}) show that these limiting behaviours on subsets $A\subset \mathbb{R}%
^{n}$ of positive $\alpha $-measure can determine global properties of these
subsets. In particular, there are bounds on the values $\alpha (A),$ $%
A\subset \mathbb{R}^{n},$ $\alpha \in \left\{ \mathcal{H}^{s}\lfloor _{A},\
P^{s}\lfloor _{A}\right\} \ $\ provided that the values $\overline{\theta }%
_{\alpha }^{s}(x)$ and $\underline{\theta }_{\alpha }^{s}(x)$ are bounded
for $x\in A$ (see \cite{MAT2}). These deep results are quite general, as
they stand for arbitrary subsets of $\mathbb{R}^{n}.$ Unfortunately, they
are useless for computing $\alpha (A)$ for general $A,$ since the
computation of the involved limits is out of reach. But if $E$ is a
self-similar set satisfying OSC and $s=\dim E$, then Mor\'{a}n proved in 
\cite{Mo1} that, 
\begin{equation}
P^{s}(E)=\sup \Big\{\left( \theta _{\alpha }^{s}(x,d)\right) ^{-1}:\text{ }%
B(x,d)\in \mathcal{B}\Big\},  \label{packing}
\end{equation}%
(recall Notation \ref{typicalball}). In that reference, analogous results
for $\alpha \in \left\{ \mathcal{H}^{s}\lfloor _{E},\mathcal{H}%
_{sph}^{s}\lfloor _{E}\right\} $ are also proved (see Sec.~\ref{section
spherical} for further details) and they were extended in \cite{LLM2} to the
characterisation of $C^{s}(E),$ namely 
\begin{equation}
C^{s}(E)=\inf \Big\{\left( \theta _{\alpha }^{s}(x,d)\right) ^{-1}:\text{ }%
x\in E\ \ \text{and \ }d>0\Big\}.  \label{centrada}
\end{equation}%
However, even with these results, the numerical computation of the values
defined in \eqref{packing} and \eqref{centrada} is still out of reach
without further restrictions on the set of balls, since the computational
time grows exponentially as their diameters decrease (see discussion in \cite%
{LLMM1}).

In \cite{LLM1, LLM3, LLM4}, it was shown that, for self-similar sets $E$
where the \textit{strong separation condition }(SSC) holds (i.e. $%
f_{i}(E)\cap f_{j}(E)=\varnothing $ for $i,j\in M,$ $i\neq j),$ the
determination of $\alpha (E),$ $\alpha \in \left\{ C^{s}\lfloor
_{E},P^{s}\lfloor _{E}\right\} ,$ starts to be computationally amenable, as
the classes of balls to be explored can be reduced to those centred at $E$
and with diameters within a known interval bounded away from zero. As the
function $\left( \theta _{\alpha }^{s}(x,d)\right) ^{-1}$ is known to be a
continuous function, both in $x$ and in $d$ (see \cite{Mat1}), the supremum
and infimum in \eqref{packing} and \eqref{centrada} then became a maximum
and a minimum, respectively. In the terminology of this article, they are $C$%
-computable. The method is not only able to render estimates of $P^{s}(E)$
and $C^{s}(E),$ but also allows an explicit construction of optimal
coverings and packings \cite{LLM1}. Moreover, explicit formulas for $\alpha
(E)$ can be found under additional, stronger forms of separation than SSC
(see \cite{Reid}).

\subsection{Results}

In this paper we make computational work with the Sierpinski gasket $S,$
where such SSC does not hold, and the results mentioned above cannot be
applied. Using the symmetries of $S,$ Theorem~\ref{Reduction1} in Sec. \ref%
{CenteredSierpinski}, narrows down the search in the class of balls given in %
\eqref{centrada} to the class of closed balls centred at points of the
cylinder set of the second generation $S_{01}$ (which amounts, from a
computational viewpoint, to $\frac{1}{9}$ of the points in the discrete
approximations of $S).$ Furthermore the existence of internal homotheties
permits the suppression of balls with a radius outside the interval $\left[ 
\frac{\sqrt{3}}{16},\frac{\sqrt{3}}{8}\right] ,$ leading to a $C$-computable
problem that is easier, in general, than the SSC case if the distance among
the $1$-cylinders is smaller than $\frac{\sqrt{3}}{16}.$ The issue of
computation of exact values of metric measures in self-similar sets with OSC
other than $S$ remains, so far, a challenge.

In Secs. \ref{subsection bounds} and \ref{subsection estimates}, we outline
the main ideas for tackling the computational task through an algorithm
aimed at approaching the value of $C^{s}(S)$ for $s=\dim S,$ with which we
obtain the following estimate 
\begin{equation*}
C^{s}(S)\sim 1.0049.
\end{equation*}%
In Theorem~\ref{theorem bounds discrete} we give, in addition to the
estimates of $C^{s}(S),$ the lower and upper bounds for $C^{s}(S)$ provided
by the algorithm at each stage $k,$ and we show that the errors tend to zero 
$(A$-computability of $C^{s}(S)).$

Notice that if we gather these results with the estimate given in \cite%
{LLMM2} for $P^{s}(S),$ we obtain full information of the total range of
values of $\theta _{\alpha }^{s}(x,d),\alpha \in \left\{ C^{s}\lfloor
_{S},P^{s}\lfloor _{S}\right\} ,$ for balls in $\mathcal{B}.$ See in Fig.\ref{optimalballs} the balls that our algorithm give as balls of maximum and
minimum density, respectively.

Finally, Sec.~\ref{section spherical} is devoted to the case of the
spherical Hausdorff measure computability problem for the Sierpinski gasket.
We have designed an algorithm that enables us to obtain an upper bound for
the value of $\mathcal{H}_{sph}^{s}(S),$ as well as to conjecture that 
\begin{equation*}
\mathcal{H}_{sph}^{s}(S)\sim 0.8616.
\end{equation*}

\section{Computability of the centred Hausdorff measure of the Sierpinski
gasket \label{CenteredSierpinski}}

\subsection{$C$-computability of $C^{s}(S)$}

The symmetry of the Sierpinski gasket can be leveraged to achieve a crucial
reduction on the set of the candidate balls to be optimal given in %
\eqref{centrada}, conducive to handling the computability problem of $%
C^{s}(S)$ with a suitable algorithm. The following two properties, valid for
a general set $E\subset\mathbb{R}^{n},$ are useful.

\begin{lemma}
\label{basiclemma} Let $E\subset\mathbb{R}^{n},$ $\alpha\in\mathcal{M}%
^{s}\lfloor_{E}$ and $h:\mathbb{R}^{n}\rightarrow\mathbb{R}^{n}$ be a
similarity with scaling factor $r_{h}.$\newline
(i) If $A\subset\mathbb{R}^{n}$ satisfies $h(A\cap E)=h(A)\cap E,$ then $%
\alpha(h(A))=r_{h}^{s}\alpha(A).$\newline
(ii) If $C\subset E$ satisfies $h(C)\subset E,$ then 
\begin{equation*}
\alpha(B(h(x),r_{h}d)\cap h(C))=r_{h}^{s}\alpha(B(x,d)\cap C).
\end{equation*}

\begin{proof}
It is well known that all metric measures scale under similarities so that,
if $\beta\in\mathcal{M}^{s},$ $h$ is a similarity with contraction ratio $%
r_{h}$ and $A\subset\mathbb{R}^{n},$ then $\beta(h(A))=r_{h}^{s}\beta(A).$
Under our hypothesis 
\begin{equation*}
\beta\lfloor_{E}(h(A))=\beta\left( h(A)\cap E\right) =\beta(h\left( A\cap
E\right) )=r_{h}^{s}\beta\lfloor_{E}(A).
\end{equation*}
This proves (i). Under the hypotheses in (ii), taking $B(x,d)\cap C$ as $A$
in (i) 
\begin{align*}
h(B(x,d)\cap C\cap E) & =h(B(x,d)\cap C)=h(B(x,d))\cap h(C) \\
& =h(B(x,d))\cap h(C)\cap E=h\left( B(x,d)\cap C)\right) \cap E,
\end{align*}
and (ii) follows.
\end{proof}
\end{lemma}

\begin{remark}
\label{mubasiclemma}If we take $E=S$ in Lemma~\ref{basiclemma}, then $\alpha$
can be replaced with $\mu,$ since $\mu$ is a multiple of the measures in $%
\mathcal{M}^{s}\lfloor_{S}.$ In general, if $E$ is a self-similar set with
OSC and $s=\dim E,$ then, $\alpha$ can be taken to be the \emph{invariant
measure of }$E$ (see the definition of invariant measure for a general
self-similar set \thinspace$E$ in \cite{HUTCH}).
\end{remark}

From now on we set $s=\dim S=\frac{\log3}{\log2}.$

The proof of Theorem \ref{Reduction1} below relies upon finding a reduced
set of balls where all the relevant values of $\theta_{\mu}^{s}(x,d)$ are
attained. A basic tool to this end is the notion of density equivalent balls.

\begin{definition}
\label{equiv} We say that $B(x,d)$ with $x\in S$ is \emph{density equivalent}
to $B(x^{\prime},d^{\prime}),$ if $x^{\prime}\in S$ and 
\begin{equation*}
\theta_{\mu}^{s}(x,d)=\theta_{\mu}^{s}(x^{\prime},d^{\prime}).
\end{equation*}
\end{definition}

The next theorem shows that $C^{s}(S)$ is $C$-computable.

\begin{theorem}
\label{Reduction1}%
\begin{equation}
C^{s}(S)=\min\left\{ \left( \theta_{\mu}^{s}(x,d)\right)^{-1}:\ \text{\ }%
x\in S_{01},\ \frac{\sqrt{3}}{16}\leq d\leq\frac{\sqrt{3}}{8}.\right\}
\label{result-Cs}
\end{equation}

\begin{proof}
We show first that 
\begin{equation}
C^{s}(S)=\inf \left\{ \left( \theta _{\mu }^{s}(x,d)\right) ^{-1}:\ \text{\ }%
x\in S_{01},\ d>0\right\}   \label{reduction1b}
\end{equation}%
Let $\pi _{i}$ be the reflection across the altitude, $h_{i},$ of $T$
(recall Notation \ref{convexhull}) through $z_{i},$ $i\in M$ (see Fig. \ref{Feasibleopenset}). Since $\pi _{1}(S_{2})=S_{0},$ $\pi _{2}(S_{1})=S_{0}$
and $\pi _{0}(S_{02})=S_{01},$ Lemma \ref{basiclemma}{\ (i)} applied to $\pi
_{i}$ and $\mu $ (see also Remark \ref{mubasiclemma}), implies, on the one
hand, that any ball centred in $S_{i},$ $i=1,2$ is density equivalent to a
ball centred in $S_{0}$ and, on the other hand, that any ball centred in $%
S_{02}$ is density equivalent to a ball centred in $S_{01}.$ So, we can
restrict our search to balls centred in $S_{00}\cup S_{01}.$

Let $x\in S_{00}$ and $d>0.$ In order to show that we can neglect $B(x,d),$
consider first the case $B(x,d)\cap S\subset S_{0}.$ Then Lemma \ref%
{basiclemma}{(ii)} implies that $B(x,d)$ is density equivalent to $%
B(f_{0}^{-1}(x),2d).$ Further, if $f_{0}^{-1}(x)\in S_{02},$ then $B(x,d)$
is density equivalent to $B(\pi_{0}(f_{0}^{-1}(x)),2d),$ centred in $S_{01}. 
$ If $f_{0}^{-1}(x)\in S_{00},$ we can repeat the argument $k$ times until $%
f_{0}^{-k}(x)\in S_{02}\cup S_{01}$ or $B(f_{0}^{-k}(x),2^{k}d)\cap(S_{1}%
\cup S_{2})\neq\emptyset.$

Consider now the case $x\in S_{00}$ and $B(x,d)\cap (S_{1}\cup S_{2})\neq
\emptyset .$ Then, $B(x,d)$ is either not optimal or it is
density-equivalent to some ball centred in $S_{01}$ because 
\begin{equation}
\left( \theta _{\mu }^{s}(x,d)\right) ^{-1}=\frac{(2d)^{s}}{\mu (B(x,d))}%
\geq \frac{(2d)^{s}}{\mu (B(\pi _{4}(x),d))}=\left( \theta _{\mu }^{s}(\pi
_{4}(x),d)\right) ^{-1},  \label{trans}
\end{equation}%
where $\pi _{4}$ is the reflection across the altitude $h_{4}$ of $T_{0}$
through $f_{2}(z_{0})$ (see Fig. \ref{Feasibleopenset}). The inequality in (%
\ref{trans}) follows from 
\begin{equation}
\mu (B(x,d))\leq \mu (B(\pi _{4}(x),d)).  \label{inequality}
\end{equation}%
In order to check \eqref{inequality}, we decompose $B(x,d)$ into the union 
\begin{equation}
B(x,d)\cap S=\left( B(x,d)\cap S_{0}\right) \cup \left( B(x,d)\cap \left(
S_{1}\cup S_{2}\right) \right)   \label{bddecomposition}
\end{equation}%
and analogously, 
\begin{equation}
B(\pi _{4}(x),d)\cap S=\left( B(\pi _{4}(x),d)\cap S_{0}\right) \cup \left(
B(\pi _{4}(x),d)\cap \left( S_{1}\cup S_{2}\right) \right) .
\label{simmetricdecomposition}
\end{equation}%
Since $\pi _{4}(S_{0})=S_{0},$ Lemma~\ref{basiclemma}{(ii)} gives 
\begin{equation*}
\mu \left( B(x,d)\cap S_{0}\right) =\mu \left( B(\pi _{4}(x),d)\cap
S_{0}\right) .
\end{equation*}%
Now, if $y\in B(x,d)\cap \left( S_{1}\cup S_{2}\right) ,$ since $S_{1}\cup
S_{2}$ is contained in the right hand half plane determined by $h_{4}$ and $x
$ belongs to the left-hand half plane, we see that $\left\vert \pi
_{4}(x)-y\right\vert \leq \left\vert x-y\right\vert $ holds, and 
\begin{equation}
B(x,d)\cap \left( S_{1}\cup S_{2}\right) \subset B(\pi _{4}(x),d)\cap \left(
S_{1}\cup S_{2}\right) .  \label{ball inclusion}
\end{equation}%
Thus, in the decompositions given in (\ref{bddecomposition}) and (\ref%
{simmetricdecomposition}) the $\mu $-measure of the first terms are equal,
whilst, by (\ref{ball inclusion}), the $\mu $-measure of the second term
cannot be smaller in (\ref{simmetricdecomposition}) than in (\ref%
{bddecomposition}), which gives (\ref{inequality}).

Once we have proved that \eqref{reduction1b} holds, we show that we can
restrict the search to balls with radii within the range $[\frac{\sqrt{3}}{16%
},\frac{\sqrt{3}}{8}].$

Let $x\in S_{01}$ and $d<\frac{\sqrt{3}}{16}.$ Clearly, $B(x,d)\cap S\subset
S_{0}\cup S_{1}.$

Suppose first that $B(x,d)\cap S\subset S_{0},$ then we can apply Lemma~\ref%
{basiclemma} {(ii)} to conclude that $B(f_{0}^{-1}(x),2d),$ centred in $S_{1}
$, is density equivalent to $B(x,d).$ Moreover, we have already seen that
such a ball is either density equivalent to a ball of equal or greater
radius and centred in $S_{01}$ or cannot be optimal. Hence, we can iterate
the argument till $d\geq\frac{\sqrt{3}}{16}$ or we find a density equivalent
ball centred in $S_{01}$ and intersecting $S_{1}.$ Observe that if we need $k
$ iterations of the argument to achieve a ball of radius $2^{k}d\geq\frac{%
\sqrt{3}}{16}$, then $2^{k-1}d<\frac{\sqrt{3}}{16}\leq2^{k}d,$ implying that 
$2^{k}d\in\lbrack\frac{\sqrt{3}}{16},\frac{\sqrt{3}}{8}].$

Now, if $B(x,d)\cap S_{1}\neq \emptyset $ (recall that we have assumed that $%
d<\frac{\sqrt{3}}{16}),$ then $x\in S_{011}$ and $B(x,d)\cap S=B(x,d)\cap
(S_{01}\cup S_{100}).$ Therefore, if we take the homothety, $h:\mathbb{R}%
^{2}\rightarrow \mathbb{R}^{2},$ of ratio $2$ and fixed point at the unique
point of $S_{0}\cap S_{1},$ $h(S_{01})=S_{0}$ and $h(S_{100})=S_{10}$ imply
that 
\begin{align*}
h\left( B(x,d)\cap S\right) & =h(B(x,d)\cap (S_{01}\cup
S_{100}))=B(h(x),2d)\cap \left( S_{0}\cup S_{10}\right) \\
& =B(h(x),2d)\cap S.
\end{align*}

Lemma~\ref{basiclemma}{(i)} shows then that $\mu(B(h(x),2d))=2^{s}%
\mu(B(x,d)),$ which implies that $B(h(x),2d)$ is density equivalent to $%
B(x,d)$. If $\frac{\sqrt{3}}{16}\leq2d,$ we have concluded because $%
2d\in\lbrack\frac{\sqrt{3}}{16},\frac{\sqrt{3}}{8}]$ and $h(x)\in S_{01}$.
Otherwise, we can repeat the argument as many times as needed until $\frac{%
\sqrt{3}}{16}\leq2d.$

Finally, if $x\in S_{01}$ and $d>\frac{\sqrt{3}}{8},$ then, since $%
f_{1}(B(x,d)\cap S)=B(f_{1}(x),\frac{d}{2})\cap S_{1}$ and $f_{1}(x)\in
S_{101}$ we get, by Lemma \ref{basiclemma} (ii), that 
\begin{align}
\mu(B(x,d)) & =\mu(B(x,d)\cap S)=2^{s}\mu(B(f_{1}(x),\frac{d}{2})\cap S_{1})
\notag \\
& \leq2^{s}\mu(B(f_{1}(x),\frac{d}{2}))=2^{s}\mu(B(\pi_{2}(f_{1}(x)),\frac {d%
}{2})),  \label{fina}
\end{align}
with $\pi_{2}(f_{1}(x))\in S_{01}.$

The proof of (\ref{result-Cs}) concludes by noticing that \eqref{fina}
implies that $(\theta_{\mu}^{s}(x,d))^{-1}\geq(\theta_{\mu}^{s}(%
\pi_{2}(f_{1}(x)),\frac{d}{2}))^{-1}$ and that we can repeat this procedure $%
k\in\mathbb{N}^{+}$ times until $\frac{\sqrt{3}}{16}\leq2^{-k}d\leq\frac{%
\sqrt{3}}{8},$ obtaining on each step a ball centred in $S_{01}$ with less
or equal inverse density.
\end{proof}
\end{theorem}

\subsection{$A$-computability of $C^{s}(S)$\label{subsection bounds}}

In this section we construct a discrete algorithm that converges to $%
C^{s}(S) $ and we provide error bounds tending to zero for its estimates,
thus showing that $C^{s}(S)$ is $A$-computable.

Following the structure of the algorithms developed in \cite{LLM4,LLM1,LLM2}%
, the construction of such an algorithm relies upon the relationship between
the centred Hausdorff measure and the inverse density function given in
Theorem~\ref{Reduction1}. With the aim of finding a computationally adequate
estimate of the minimum value given in \eqref{result-Cs}, a discrete
approximation of both the Sierpinski gasket and the invariant measure is
proposed.

Firstly, it is well-known that, for any non-empty compact subset $A\subset%
\mathbb{R}^{2},$ $S$ can be built with an arbitrary level of detail by
increasing the iterations $k$ in $F^{k}(A),$ where $F^{k}=F\circ F...\circ F$
is the $k$-th iterate of the Hutchinson operator $F$ (see \eqref{hutchin}).
This is because $\lim_{k\rightarrow\infty}F^{k}(A)=S$ with respect to the
Hausdorff metric, given that $S$ is the attractor of $\Psi$ under the
contracting operator $F$ (cf. \cite{HUTCH}). Furthermore, if $A\subset S,$
then $F^{k}(A)\subset S$ for any $k\in\mathbb{N^{+}}.$ We use these facts in
the design of our algorithm where we take as initial compact set $%
A_{1}:=\{z_{0},z_{1},z_{2}\}$ (recall that $z_{i}\in M$ are the fixed points
of the similitudes in $\Psi)$ and obtain the set 
\begin{equation}
A_{k}:=F^{k-1}(A_{1})\subset S,\text{ }k\geq2  \label{Ak}
\end{equation}
as a discrete approximation of $S$ at iteration $k$ of our algorithm.

If we take $\alpha =\mu _{1}:=\frac{1}{3}(\delta _{z_{0}}+\delta
_{z_{1}}+\delta _{z_{2}})$ as the initial measure in (\ref{invariant}),
where $\delta _{x}$ is the Dirac probability measure at $x,$ then 
\begin{equation}
\mu _{k}:=\mathbf{M}^{k-1}(\mu _{1})=\frac{1}{3^{k-1}}\sum_{i\in M^{k-1}}\mu
_{1}\circ f_{i}^{-1}=\frac{1}{3^{k}}\sum_{i\in M^{k-1}}\left( \delta
_{f_{i}(z_{0})}+\delta _{f_{i}(z_{1})}+\delta _{f_{i}(z_{2})}\right) 
\label{mu(k)}
\end{equation}%
is a probability measure supported on $A_{k}\subset S$ and $\mu _{k}%
\xrightarrow[k \to \infty]{w}\mu .$

The discrete measure $\mu_{k}$ is the approximation of the invariant measure 
$\mu$ that our algorithms take at iteration $k.$

For $i\in M^{k-1},$ $\delta _{f_{i_{1}...i_{k-2}}(z_{i_{k-1}})}$ is an atom
of $\mu _{k}$ and a summand in the right-hand term in (\ref{mu(k)}). Since $%
f_{i}(z_{j})=f_{j}(z_{i}),$ for $i,$ $j\in M,$ with $i\neq j,$ all the
points in $A_{2}-A_{1}$ have two codes in $M^{2}.$ From this, it easily
follows that all the points in $A_{k}-A_{1}$ also have two codes in $M^{k},$ 
$k\geq 2,$ (see \cite{LLMM2}, Sec. 3), and therefore we can write (\ref%
{mu(k)}) as%
\begin{equation}
\mu _{k}=\frac{1}{3^{k}}(\delta _{z_{0}}+\delta _{z_{1}}+\delta _{z_{2}})+%
\frac{2}{3^{k}}\sum_{x\in A_{k}\setminus A_{1}}\delta _{x}.  \label{mu(k)b}
\end{equation}%
The algorithm outlined in Sec.~\ref{subsection estimates} works with the
sets $A_{k}$ $\subset S$ defined in (\ref{Ak}) and with the measures $\mu
_{k}$ defined in (\ref{mu(k)b}) as approximations of the Sierpinski gasket $%
S $ and of the invariant measure $\mu $ at iteration $k,$ respectively.

Addressing the issue of establishing error bounds in the estimates of the $%
\mu$-measure of balls through the approximating $\mu_{k}$-measures, the
following lemma allows us to make the comparison of the $\mu$ measure and
the $\mu_{k}$-measure in two relevant cases.

\begin{lemma}
\label{union of cylinders}$\ \ $\newline
(i) Let $\left\{ S_{i}:i\in I\subset M^{k}\right\} ,$ $k\in\mathbb{N}^{+},$
be a collection of $k$-cylinder sets. Then%
\begin{equation*}
\mu\left({\displaystyle\bigcup\limits_{i\in I}}S_{i}\right)
\leq\mu_{k}\left( {\displaystyle\bigcup\limits_{i\in I}}S_{i}\right)
\end{equation*}
(ii) Let $A\subset S,$ $k\in\mathbb{N}^{+},$ and $I:=\left\{i\in
M^{k}:S_{i}\cap A\right\} \neq\varnothing.$ Then

\begin{equation*}
\mu_{k}(A)\leq\mu\left({\displaystyle\bigcup\limits_{i\in I}}S_{i}\right)
\end{equation*}

\begin{proof}
For any $i\in I$ we know that $\mu (S_{i})=$ $3^{-k},$ so%
\begin{equation}
\mu \left( {\displaystyle\bigcup\limits_{i\in I}}S_{i}\right) =\sum_{i\in
I}\mu \left( S_{i}\right) =3^{-k}\#(I)  \label{mucilindros}
\end{equation}%
where $\#(I)$ denotes the cardinality of $I.$ We know that each cylinder $%
S_{i}$ contains a unique point, $\hat{x}$ in $A_{k}$ (namely $\hat{x}:=\pi
(i_{1}i_{2}...i_{k-1}i_{k}i_{k}i_{k}...))=f_{i_{1}i_{2}...i_{k-1}}(z_{i_{k}})),
$ and that each point of $({\displaystyle\bigcup\limits_{i\in I}}S_{i})\cap
A_{k}$ either belongs to a unique $k$-cylinder in $\left\{ S_{i}:i\in
I\right\} $ (let us write $A_{k,1}(I)$ for such subset of $A_{k}),$ or
belongs to the set $A_{k,2}(I),$ or subset of points of $A_{k}$ that lay on
two $k$-cylinders of $\left\{ S_{i}:i\in I\right\} .$ Hence $A_{k}\cap ({%
\displaystyle\bigcup\limits_{i\in I}}S_{i})=A_{k,1}(I)\cup A_{k,2}(I)$ with $%
A_{k,1}(I)\cap A_{k,2}(I)=\varnothing .$ Then, since $A_{1}\cap ({%
\displaystyle\bigcup\limits_{i\in I}}S_{i})\subset A_{k,1}(I)\cap ({%
\displaystyle\bigcup\limits_{i\in I}}S_{i})$ and, in consequence $A_{1}\cap
A_{k,2}(I)=\varnothing ,$ we have 
\begin{align}
\mu _{k}\left( {\displaystyle\bigcup\limits_{i\in I}}S_{i}\right) & =\mu
_{k}\left( \left( {\displaystyle\bigcup\limits_{i\in I}}S_{i}\right) \cap
A_{k}\right)   \notag \\
& =\mu _{k}\left( A_{k,1}(I)\cap A_{1}\right) +\mu _{k}\left( A_{k,1}(I)\cap
(A_{k}-A_{1})\right) +\mu _{k}\left( A_{k,2}(I)\cap (A_{k}-A_{1})\right)  
\notag \\
& =\frac{1}{3^{k}}\left\{ \#\left( A_{k,1}(I)\cap A_{1}\right) +2\#\left(
A_{k,1}(I)\cap (A_{k}-A_{1})\right) +2\#\left( A_{k,2}(I)\cap
(A_{k}-A_{1}\right) )\right\}   \notag \\
& \geq \frac{1}{3^{k}}\left\{ \#\left( A_{k,1}(I)\right) +2\#\left(
A_{k,2}(I)\cap (A_{k}-A_{1}\right) )\right\} =\frac{1}{3^{k}}\left\{
\#\left( A_{k,1}(I)\right) +2\#\left( A_{k,2}(I)\right) \right\} .
\label{partition1}
\end{align}%
Let $J:=\{i\in M^{k}:i=i_{1}i_{1}...i_{1},$ $i_{1}\in M\},$ $k>1$ and $i\in
M^{k}\backslash J.$ For each cylinder $S_{i},$ there is $i^{\ast }\in
M^{k}\backslash J$ such that $S_{i}\cap S_{i^{\ast }}=S_{i}\cap S_{i^{\ast
}}\cap A_{k}$ consists of a unique point. Consider the partition of $I,$ $%
I=I_{0}\cup I_{1}\cup I_{2}$ where $I_{0}:=J\cap I,$ $I_{1}:=\{i\in I:\ \
i^{\ast }\notin I\}$ and $I_{2}:=\{i\in I:\ \ i^{\ast }\in I\}.$ Since there
is a bijective mapping between $A_{k,1}(I)$ and $I_{0}\cup I_{1},$ an
injective mapping between $A_{k,2}(I)$ and $I_{2},$ and each cylinder $S_{i},
$ $i\in M^{k}$ contains a unique point in $A_{k,2}(I),$ we get that $%
\#(A_{k,1}(I))=\#(I_{0})+(\#(I_{1}))$ and $2\#(A_{k,2}(I))=(\#(I_{2})).$ This,
together with (\ref{partition1}) and (\ref{mucilindros}) gives 
\begin{equation*}
\mu _{k}\left( {\displaystyle\bigcup\limits_{i\in I}}S_{i}\right) \geq \frac{%
1}{3^{k}}\#\left( I\right) =\mu \left( {\displaystyle\bigcup\limits_{i\in I}}%
S_{i}\right) .
\end{equation*}

Observe that the equality holds if and only if the set in the other term
affected by a coefficient 2, $A_{k,1}(I)\cap(A_{k}-A_{1}),$ is empty.

In order to prove (ii), note that $A_{k,1}(I)\cap (A_{k}-A_{1})=\varnothing ,
$ since if some cylinder in $\left\{ S_{i}:i\in I\right\} $ intersects the
set $A\cap (A_{k}-A_{1})$ at a point $x,$ then there is another cylinder to
which $x$ also belongs and this cylinder also belongs in turn to the
collection of cylinders $\left\{ S_{i}:i\in I\right\} .$ Thus, using the
above decomposition of $\mu _{k}\left( \left( {\displaystyle%
\bigcup\limits_{i\in I}}S_{i}\right) \cap A_{k}\right) $ and by the final
observation, 
\begin{align*}
\mu _{k}(A)& =\mu _{k}(A\cap A_{k})\leq \mu _{k}\left( ({\displaystyle%
\bigcup\limits_{i\in I}}S_{i})\cap A_{k}\right)  \\
& =\mu _{k}\left( A_{k,1}(I)\cap A_{1}\right) +\mu _{k}\left( A_{k,2}(I)\cap
(A_{k}-A_{1})\right)  \\
& =\frac{1}{3^{k}}\#\left( I\right) =\mu ({\displaystyle\bigcup\limits_{i\in
I}}S_{i}).
\end{align*}
\end{proof}
\end{lemma}

Theorem \ref{theorem bounds discrete} below establishes the $A$%
-computability of $C^{s}(S).$ It is the discrete version of Theorem \ref%
{Reduction1}, suitable for our computational purposes as it gives the
estimates and the lower and upper bounds of $C^{s}(S)$ at each iteration $k,$
namely, $C_{k},$ $C_{k}^{\inf}$ and $C_{k}^{\sup}$ respectively. We first
prove the following lemma.

\begin{lemma}
\label{lemma cotas}Let $k>0,$ $x\in 
\mathbb{R}
^{2}$ and $d>2^{-k}.$ Then,$\,\ $\newline
(i) $\mu (B(x,d))\geq \mu _{k}(B(x,d-2^{-k}))$\newline
(ii) If $S\nsubseteq B(x,d)$ and $B(x,d)\cap A_{k}\neq \varnothing ,$ then
there is $y_{k}\in A_{k}$ satisfying $\mu (B(x,d))\leq \mu _{k}(B(x,d_{k})),$
where $d_{k}=\left\vert y_{k}-x\right\vert $ and $d-2^{-k}\leq d_{k}\leq
d+2^{-k}.$

\begin{proof}
Let $k>0,$ $x\in 
\mathbb{R}
^{2}$ and $d>2^{-k}.$ \newline
(i) Let 
\begin{equation*}
J_{k}:=\{i\in M^{k}:S_{i}\subset B(x,d)\}
\end{equation*}%
and 
\begin{equation*}
H_{k}:=\{i\in M^{k}:B(x,d-2^{-k})\cap S_{i}\neq \emptyset \}.
\end{equation*}%
Clearly, $H_{k}\subset J_{k}$ holds. Then, using Lemma \ref{union of
cylinders} (ii) 
\begin{equation*}
\mu _{k}(B(x,d-2^{-k}))\leq \sum_{i\in H_{k}}\mu (S_{i})\leq \sum_{i\in
J_{k}}\mu (S_{i})\leq \mu (B(x,d)).
\end{equation*}%
(ii) We prove first that, if $S\nsubseteq B(x,d)$ and $B(x,d)\cap A_{k}\neq
\varnothing ,$ the set 
\begin{equation*}
G_{k}:=\{i\in M^{k}:\partial B(x,d)\cap S_{i}\neq \emptyset \}
\end{equation*}%
is non empty, where $\partial B(x,d)$ is the border of $B(x,d).$ Let $U(x,d)$
denote the open ball centred at $x$ and with radius $d.$ If $%
G_{k}=\varnothing ,$ then 
\begin{equation}
{\displaystyle\bigcup\limits_{i\in M^{k}}}S_{i}\subset U(x,d)\cup \left(
B(x,d)\right) ^{c}.  \label{disconnected}
\end{equation}%
Since $B(x,d)\cap A_{k}\neq \varnothing $ and we may assume that $\partial
B(x,d)\cap A_{k}=\varnothing $ or $G_{k}$ would be trivially nonempty, we
see that $U(x,d)\cap (\cup _{i\in M^{k}}S_{i})\neq \varnothing ,$ and $%
S\nsubseteq B(x,d)$ implies that the set $\left( B(x,d)\right) ^{c}\cap
(\cup _{i\in M^{k}}S_{i})$ is nonempty, so (\ref{disconnected}) contradicts
that for any $k>0,$ $\cup _{i\in M^{k}}S_{i}$ is a connected set and our
claim that $G_{k}\neq \varnothing $ is proved. \newline
Let 
\begin{equation*}
I_{k}:=\{i\in M^{k}:B(x,d)\cap S_{i}\neq \emptyset \}
\end{equation*}%
and%
\begin{equation*}
d_{k}:=\max \{\left\vert y-x\right\vert :\text{ }y\in A_{k}\cap \left( \cup
_{i\in I_{k}}S_{i}\right) \}.
\end{equation*}%
Let $z\in A_{k}$ satisfying $d_{k}=\left\vert z-x\right\vert .$ Then, $%
d_{k}\leq d+2^{-k},$ since $\cup _{i\in I_{k}}S_{i}\subset B(x,d+2^{-k})$
and, using part (i) of Lemma \ref{union of cylinders}, we have 
\begin{align*}
\mu (B(x,d))& \leq \mu (\cup _{i\in I_{k}}S_{i})\leq \mu _{k}(\cup _{i\in
I_{k}}S_{i}) \\
& =\mu _{k}(\left( \cup _{i\in I_{k}}S_{i}\right) \cap A_{k})\leq \mu
_{k}\left( B(x,d_{k})\right) .
\end{align*}%
Finally, using that $\partial B(x,d)\cap (\cup _{i\in M^{k}}S_{i})\neq
\varnothing $ we see that 
\begin{equation*}
d_{k}\geq \max \{\left\vert y-x\right\vert :\text{ }y\in A_{k}\cap \left(
\cup _{i\in G_{k}}S_{i}\right) \}\geq d-2^{-k}.
\end{equation*}
\end{proof}
\end{lemma}

\begin{theorem}
\label{theorem bounds discrete} Let 
\begin{equation}
C_{k}=\min \left\{ (\theta _{\mu _{k}}^{s}(x,d))^{-1}:\text{ }x\in A_{k}\cap
S_{01}\text{ and }d=\left\vert y-x\right\vert \text{ with }y\in A_{k}\text{
and \ }\frac{\sqrt{3}}{16}\leq d\leq \frac{\sqrt{3}}{8}+2^{1-k}\right\} .
\label{Ck}
\end{equation}%
Then, for every $k\geq 4,$ 
\begin{equation}
C_{k}^{\inf }\leq C^{s}(S)\leq C_{k}^{\sup }  \label{Csbounds}
\end{equation}%
holds, where 
\begin{equation}
C_{k}^{\inf }=K_{k}C_{k},\text{ \ }K_{k}=\left( 1+\frac{2^{5-k}}{\sqrt{3}}%
\right) ^{-s},\text{ \ \ }C_{k}^{\sup }=\frac{(2d_{k})^{s}}{\mu
_{k}(B(x_{k},d_{k}-2^{-k}))}  \label{Cs bounds and k}
\end{equation}%
and $B(x_{k},d_{k})$ is a ball minimising (\ref{Ck}).

\begin{proof}
Let $k\geq 4,$ and let $B(x_{k},d_{k})$ be a ball minimising (\ref{Ck}). We
know that $C_{k}=\left( \theta _{\mu _{k}}^{s}(x_{k},d_{k}\right) )^{-1},$
and $\frac{\sqrt{3}}{16}\leq d_{k}\leq \frac{\sqrt{3}}{8}+2^{1-k}.$ Using %
\eqref{centrada} and Lemma \ref{lemma cotas} (i), we get 
\begin{equation*}
C^{s}(S)\leq \frac{(2d_{k})^{s}}{\mu (B(x_{k},d_{k}))}\leq \frac{(2d_{k})^{s}%
}{\mu _{k}(B(x_{k},d_{k}-2^{-k}))}=C_{k}^{\sup },
\end{equation*}%
where $C_{k}^{\sup }$ is well defined since $d_{k}-2^{-k}\geq \frac{\sqrt{3}%
}{16}-2^{-k}\geq 0$ and $\mu _{k}(B(x_{k},d_{k}-2^{-k}))>0,$ as $x_{k}\in
B(x_{k},d_{k}-2^{-k})\cap A_{k}.$\newline
In order to prove the inequality $C_{k}^{\inf }\leq C^{s}(S),$ let $(x,d)\in
S\times \lbrack \frac{\sqrt{3}}{16},\frac{\sqrt{3}}{8}]$ be such that $%
C^{s}(S)=\frac{(2d)^{s}}{\mu (B(x,d))}.$ Let $i\in M^{k}$ be such that $x\in
S_{i},$ and let $y_{k}=S_{i}\cap A_{k}.$ Since $B(x,d)\subset
B(y_{k},d+2^{-k})$ we get 
\begin{equation}
\mu (B(x,d))\leq \mu (B(y_{k},d+2^{-k})).  \label{ineq}
\end{equation}%
Taking the ball $B(y_{k},d+2^{-k})$ as $B(x,d)$ in Lemma \ref{lemma cotas}
(ii) we can get a point $z_{k}\in A_{k}$ satisfying that 
\begin{equation}
\mu (B(y_{k},d+2^{-k}))\leq \mu _{k}(B(y_{k},d_{k}^{\ast })),  \label{ineq1}
\end{equation}%
where $d_{k}^{\ast }:=\left\vert z_{k}-y_{k}\right\vert $ and $d_{k}^{\ast
}\in \lbrack d,d+2^{1-k}]\subset \lbrack \frac{\sqrt{3}}{16},\frac{\sqrt{3}}{%
8}+2^{1-k}].$ Then using (\ref{ineq}), (\ref{ineq1}), (\ref{Ck}), that $%
d_{k}^{\ast }\leq d+2^{1-k},$ and that $d\geq \frac{\sqrt{3}}{16}$ we get 
\begin{gather*}
C^{s}(S)=\frac{(2d)^{s}}{\mu (B(x,d))}\geq \frac{(2d)^{s}}{\mu
_{k}(B(y_{k},d_{k}^{\ast }))}=\Big(\frac{d_{k}^{\ast }}{d}\Big)^{-s}\left(
\theta _{\mu _{k}}^{s}(y_{k},d_{k}^{\ast })\right) ^{-1} \\
\geq \Big(\frac{d_{k}^{\ast }}{d}\Big)^{-s}C_{k}\geq \Big(\frac{d+2^{1-k}}{d}%
\Big)^{-s}C_{k}\geq \Big(1+\frac{2^{5-k}}{\sqrt{3}}\Big)^{-s}C_{k}=C_{k}^{%
\inf }.
\end{gather*}
\end{proof}

\begin{remark}
Notice that (\ref{Csbounds}) and (\ref{Cs bounds and k}) give the
convergence of $C_{k}$ to $C^{s}(S),$ as $k$ tends to infinity, and
therefore the convergence of the algorithm.
\end{remark}
\end{theorem}

\begin{remark}
\label{estimate center}Theorem \ref{theorem bounds discrete} gives error
bounds in the estimation of $C^{s}(S)$ in terms of the density of the ball
selected by the algorithm as optimal, but it does not give any error bound
on the estimation of the centre and of the radius of the optimal ball for
the measure $\mu .$ Such estimations remain so far an open problem.
\end{remark}

\subsection{Algorithm and numerical results\label{subsection estimates}}

The numerical results of this subsection are obtained through an algorithm
that computes $C_{k}$ and the bounds of $C^{s}(S),$ $C_{k}^{\inf}$ and $%
C_{k}^{\sup},$ given in Theorem \ref{theorem bounds discrete}.

\subsubsection{The algorithm}

The structure of this new algorithm is akin to the one presented in \cite%
{LLMM2} for $P^{s}(S),$ so we will simply describe it in general terms,
making a comparison with its equivalent for the packing measure and
referring the interested reader to \cite{LLMM2} for further details.

The analogous result to Theorem~\ref{theorem bounds discrete} for $P^{s}(S)$
is included in the remark below (Theorem 7 in \cite{LLMM2}). The proof of
Theorem~\ref{theorem bounds discrete} using part (ii) of Lemma \ref{lemma
cotas} introduces a novel approach that can be generalised to improve the
values of $d_{0},$ $K_{k}^{P},$ and the restriction $k\geq 6$ in the remark
below, which could be replaced with $\widetilde{d}_{0}:=\frac{\sqrt{3}}{16}%
-2^{1-k},$ $\widetilde{K}_{k}^{P}:=\left( 1-\frac{2^{5-k}}{\sqrt{3}}\right)
^{-s},$ and $k\geq 4$ respectively. Remark \ref{theorem packing} will ease
the understanding of the changes required for the adaptation of the
algorithm for the estimation of $P^{s}(S)$ to the estimation of $C^{s}(S).$

\begin{remark}
\label{theorem packing}For every $k\geq 6,$ 
\begin{equation*}
P_{k}^{\inf }\leq P^{s}(S)\leq P_{k}^{\sup }
\end{equation*}%
where 
\begin{equation}
P_{k}:=\max \left\{ \left( \mathring{\theta}_{\mu _{k}}^{s}(x,d)\right)
^{-1}:\text{ }x\in A_{k}\cap S_{01},\text{ }d=\left\vert y-x\right\vert 
\text{ with }y\in A_{k}\backslash S_{2}\text{ and }d_{0}\leq d\leq
d_{x}\right\} ,  \label{Pk}
\end{equation}%
$d_{0}:=\frac{\sqrt{3}}{16}-2^{2-k},$ $d_{x}:=\max \{\left\vert
y-x\right\vert :$ $y\in \partial \mathcal{R\}},$ $\mathcal{R}$ is a feasible
open set for $S,$ $\mathring{\theta}_{\mu _{k}}^{s}(x,d)$ is the $\mu _{k}$%
-density of the open ball $U(x,d),$ 
\begin{equation}
P_{k}^{\inf }:=\frac{(2d_{k})^{s}}{\mu _{k}(U(x_{k},d_{k}+2^{-k}))},\text{ \ 
}K_{k}^{P}:=\left( 1-\frac{2^{6-k}}{\sqrt{3}}\right) ^{-s},\text{ }%
P_{k}^{\sup }:=K_{k}^{P}P_{k}.  \label{bounds packing}
\end{equation}%
and $U(x_{k},d_{k})$ is an open ball that maximises (\ref{Pk}).
\end{remark}

Let us briefly list the changes required to adapt the algorithm described in 
\cite{LLMM2} and based on the result of Remark \ref{theorem packing} to our
case, avoiding unnecessary duplication.

\begin{enumerate}
\item \textbf{Replace maximums with minimums.} Recall that for any $k\geq4,$
the aim of the algorithm for the estimation of $C^{s}(S)$ is to find a ball, 
$B(x,d)$ of minimal inverse $\mu_{k}$-density where $d=\left\vert
x-y\right\vert ,$ $x\in A_{k}\cap S_{01},$ and $y\in A_{k}.$ For the
estimation of $P^{s}(S)$ the goal was to maximise the inverse density, so
this has to be changed accordingly.

\item \textbf{New bounds.} The definitions of $C_{k}^{\inf},$ $C_{k}^{\sup} $
and $K_{k}$ given in \eqref{Cs bounds and k}, are analogous to the ones of $%
P_{k}^{\inf},$ $P_{k}^{\sup}$ and $K_{k}^{P}$ given in (\ref{bounds packing}%
). Thus, to rewrite the algorithm the roles of $P_{k}^{\inf},$ $K_{k}^{P}$
and $P_{k}^{\sup}$ have to be switched with those of $C_{k}^{\inf},$ $K_{k}$
and $C_{k}^{\sup},$ respectively.

\item \textbf{Range of radii.} Computing the terms $\mu_{k}(B(x,d))$ and $%
(2d)^{s}$ needed to find the minimum value given in \eqref{Ck}, requires the
calculation of the distances from each $x\in A_{k}\cap S_{01}$ to all the
points in $A_{k}$ and the selection of those within the allowed range of
radii given in (\ref{Ck}). The $\mu_{k}$-measure of the resulting balls is
then obtained by arranging in a list the sequence of feasible distances in
increasing order so that the position of a distance $d$ in the list,
together with the distances that are equal to $d,$ provides the number of
points of $A_{k}$ within $B(x,d),$ and hence it gives $\mu_{k}(B(x,d)).$
Since the constrains on the radii of these candidates to optimal balls
differ from one case to another, it should be adapted by replacing $d_{0}:=%
\frac{\sqrt{3}}{16}-2^{2-k}$ with $\frac{\sqrt{3}}{16}$ and $d_{x}$ with $%
\frac{\sqrt{3}}{8}+2^{1-k}.$

\item \textbf{Closed balls.} The bounds in Theorem \ref{theorem bounds
discrete}, and in particular the inequality (\ref{ineq1}), require that the
balls considered by the algorithm be closed, instead of open like those used
in the computation of $P^{s}(S).$ Notice that since $\mu (\partial B(x,d))=0$
(see \cite{Mat1}), using open or closed balls in \eqref{result-Cs} does not
make a difference, and therefore neither does it make a difference in %
\eqref{Ck} or \eqref{Pk} if $k$ is large enough. Only for small $k$ do the
results vary substantially. The use of closed balls implies that, in the
number of points that contribute to the $\mu _{k}$-measure of the ball $%
B(x,d),$ we have to consider the number of points, $t_{x},$ in $U(x,d),$ and
the number of points, $T_{x},$ in $\partial B(x,d).$\ Therefore we have to
replace $\mu _{k}(U(x,d))=\frac{2}{3^{k}}t_{x}$ in \cite{LLMM2} with $\mu
_{k}(B(x,d))=\frac{2}{3^{k}}(t_{x}+T_{x}).$
\end{enumerate}

\subsubsection{Numerical results}

\begin{figure}[H]
\begin{subfigure}[b]{0.5\textwidth}
\centering
\includegraphics[width=\textwidth]{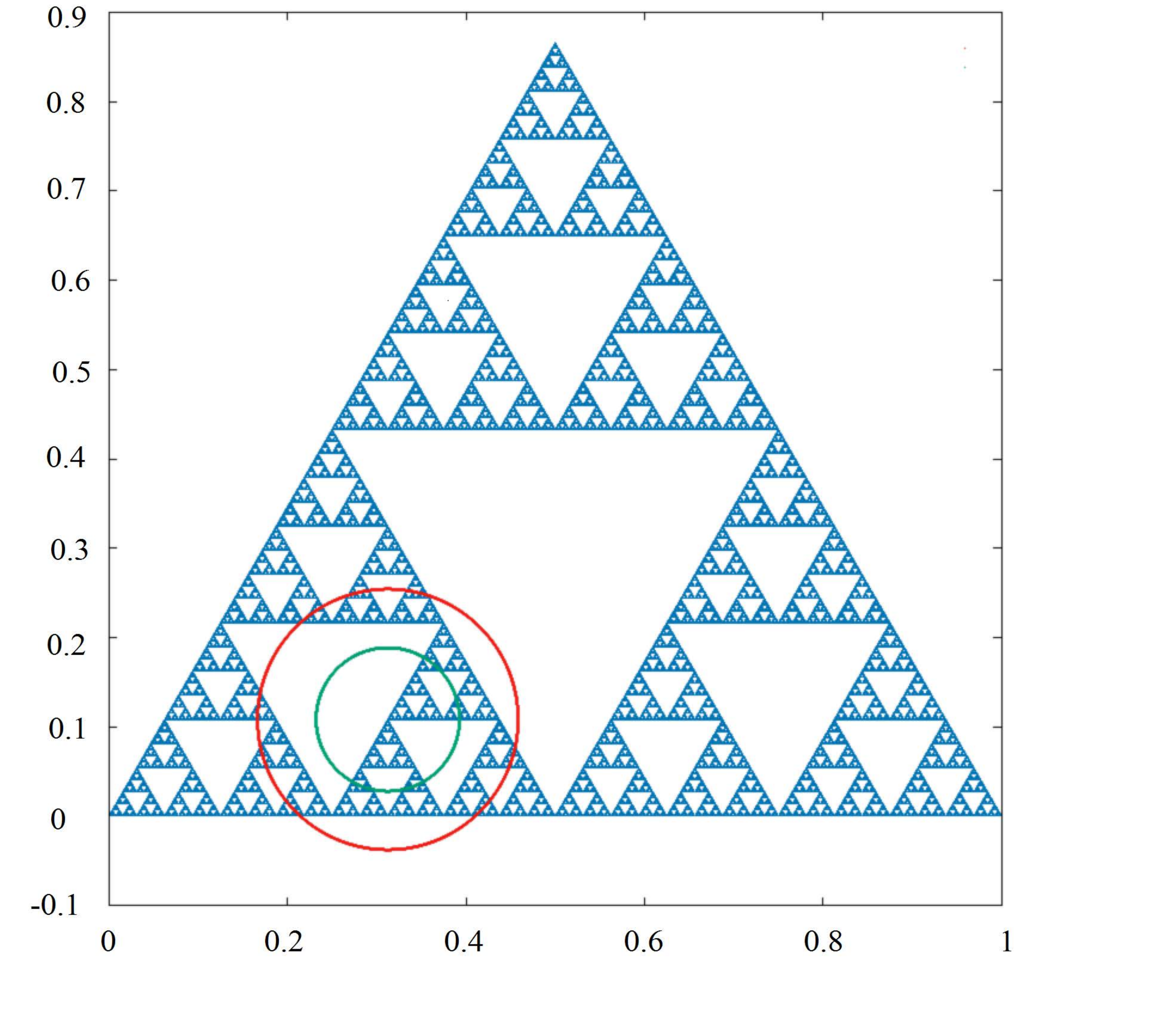}
\end{subfigure}
\hfill%
\begin{subfigure}[b]{0.5\textwidth}
\centering
\includegraphics[width=\textwidth]{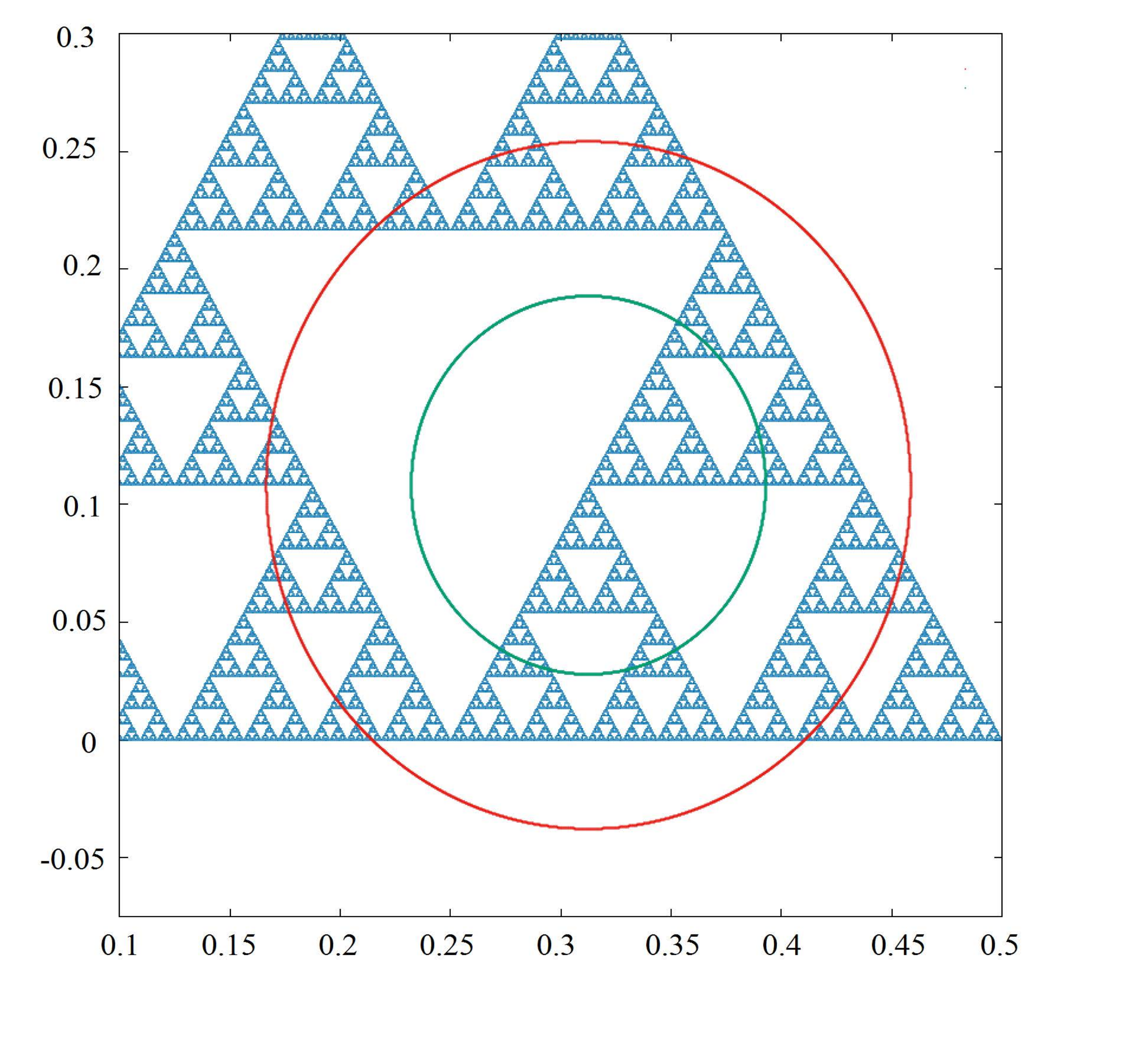}
\end{subfigure}
\caption{Optimal balls for the packing and centred measures of $S$. \newline
Balls of minimum (in green) and maximum (in red) $\protect\mu _{14}$-density.
\newline
}
\label{optimalballs}
\end{figure}

Table~\ref{numresults} shows the algorithm's output from the fifth to the
fourteenth iteration. The output of the algorithm for $k=14$ together with
Theorem~\ref{theorem bounds discrete} gives the estimate $C_{14}=1.004903$
of $C^{s}(S),$ and a $100\%$ confidence interval for $C^{s}(S),$ $%
[C_{14}^{\inf },C_{14}^{\sup }]=[1.003109,1.005611],$ with a length of less
than $0.002502.$ Moreover, for every iteration $k$ the selected ball, $%
B(x_{k},d_{k}),$ has the same centre $x_{k}:=f_{010}(z_{2})=(\frac{5}{16},%
\frac{\sqrt{3}}{16})$ and its radius, $d_{k},$ varies slightly for any $k$
in the range $k=11,...,14.$ The stability observed indicates that the
inverse $\mu _{k}$-density of the ball $B(x_{k},0.146)$ (the red colour ball
in Fig.~\ref{optimalballs}) results in a good approximation to the minimum
value in (\ref{result-Cs}).

Regarding the stability of the error bounds for $C^{s}(S),$ Table~\ref{numresults} shows that there are two fixed decimal places in the values of
the last three iterations of $C_{k}^{\sup }$, and a slightly slower
stabilisation of $C_{k}^{\inf }.$ This is mainly due to the slow convergence
to one of the terms $K_{k},$ which was also the reason why, in the algorithm
for $P^{s}(S)$, the convergence of $P_{k}^{\sup }$ was slower than that of $%
P_{k}^{\inf }$ (see Sec. 4.3 in \cite{LLMM2}).

\begin{table}[H]
\centering%
\begin{equation*}
\begin{tabular}{|c|c|c|c|c|}
\hline
$k$ & $d_{k}$ & $C_{k}^{\inf}$ & $C_{k}$ & $C_{k}^{\sup}$ \\ \hline
5 & 0.125 & 0.409736 & 0.843750 & 2.700000 \\ \hline
\multicolumn{1}{|l|}{6} & \multicolumn{1}{|l|}{0.143205} & 
\multicolumn{1}{|l|}{0.622414} & \multicolumn{1}{|l|}{0.930364} & 
\multicolumn{1}{|l|}{1.255991} \\ \hline
\multicolumn{1}{|l|}{7} & \multicolumn{1}{|l|}{0.143205} & 
\multicolumn{1}{|l|}{0.790389} & \multicolumn{1}{|l|}{0.978694} & 
\multicolumn{1}{|l|}{1.141810} \\ \hline
\multicolumn{1}{|l|}{8} & \multicolumn{1}{|l|}{0.144690} & 
\multicolumn{1}{|l|}{0.894667} & \multicolumn{1}{|l|}{0.999143} & 
\multicolumn{1}{|l|}{1.068851} \\ \hline
\multicolumn{1}{|l|}{9} & \multicolumn{1}{|l|}{0.144690} & 
\multicolumn{1}{|l|}{0.945925} & \multicolumn{1}{|l|}{1.000593} & 
\multicolumn{1}{|l|}{1.035149} \\ \hline
\multicolumn{1}{|l|}{10} & \multicolumn{1}{|l|}{0.147354} & 
\multicolumn{1}{|l|}{0.975686} & \multicolumn{1}{|l|}{1.003735} & 
\multicolumn{1}{|l|}{1.016677} \\ \hline
\multicolumn{1}{|l|}{11} & \multicolumn{1}{|l|}{0.145596} & 
\multicolumn{1}{|l|}{0.990358} & \multicolumn{1}{|l|}{1.004556} & 
\multicolumn{1}{|l|}{1.011856} \\ \hline
\multicolumn{1}{|l|}{12} & \multicolumn{1}{|l|}{0.145834} & 
\multicolumn{1}{|l|}{0.997550} & \multicolumn{1}{|l|}{1.004691} & 
\multicolumn{1}{|l|}{1.007754} \\ \hline
\multicolumn{1}{|l|}{13} & \multicolumn{1}{|l|}{0.145957} & 
\multicolumn{1}{|l|}{1.001285} & \multicolumn{1}{|l|}{1.004867} & 
\multicolumn{1}{|l|}{1.006332} \\ \hline
\multicolumn{1}{|l|}{14} & \multicolumn{1}{|l|}{0.145957} & 
\multicolumn{1}{|l|}{1.003109} & \multicolumn{1}{|l|}{1.004903} & 
\multicolumn{1}{|l|}{1.005611} \\ \hline
\end{tabular}%
\end{equation*}%
\caption{Estimates of the centred measure of $S.$\newline
Algorithm outputs rounded to six decimal places: to the smallest value for $%
C_{k}^{\inf }$, to the largest for $C_{k}^{\sup }$ and to the nearest for $%
C_{k}$ and $d_{k}.$ }
\label{numresults}
\end{table}
Since the values of $C_{k}$ change at a slower rate from $k=10$ on, and the
lower and upper bounds of $C_{k}$ are arbitrarily close to $C_{k}$ for large
enough $k,$ we can conjecture that $C^{s}(S)\sim C_{14}\sim 1.0049.$ More
precise estimates for $C^{s}(S)$ would require a refinement of the lower
bound for $C^{s}(S)$ or a significant increment of the largest value taken
by $k,$ now fixed at $k_{\max }:=14.$

If we gather these results with the estimate given in \cite{LLMM2} for $%
P^{s}(S),$ we obtain a quite complete information of the total range of
values of $\theta_{\alpha}^{s}(x,d),$ $\alpha\in\left\{
C^{s}\lfloor_{S},P^{s}\lfloor_{S}\right\}.$

In \cite{LLMM2} what is obtained is the estimate $P_{15}=1.6683$ for $%
P^{s}(S)$ at iteration $k=15,$ and balls $B(z_{k},d_{k})$ maximising (\ref%
{Pk}) for $k\in \{6,...,15\}$ were found. For $k=14,15$ the centres of these
balls are both the same $z_{k}:=(0.5,0),$ and their radii are $d_{k}\sim
0.1605.$ Notice that the ball $B(z_{k},0.1605)$ is $\mu $-density equivalent
to a ball centred at $f_{010}(z_{2}),$ which is precisely the centre of the
ball which minimises (\ref{Ck}), and with radius $d\sim 0.08.$ This ball is
outlined in green in Fig.~\ref{optimalballs}, together with the red colour
ball whose inverse density gives the estimate of $C^{s}(S)$ at iteration $%
k=14.$

\begin{figure}[H]
\centering
\includegraphics[width=0.6\textwidth]{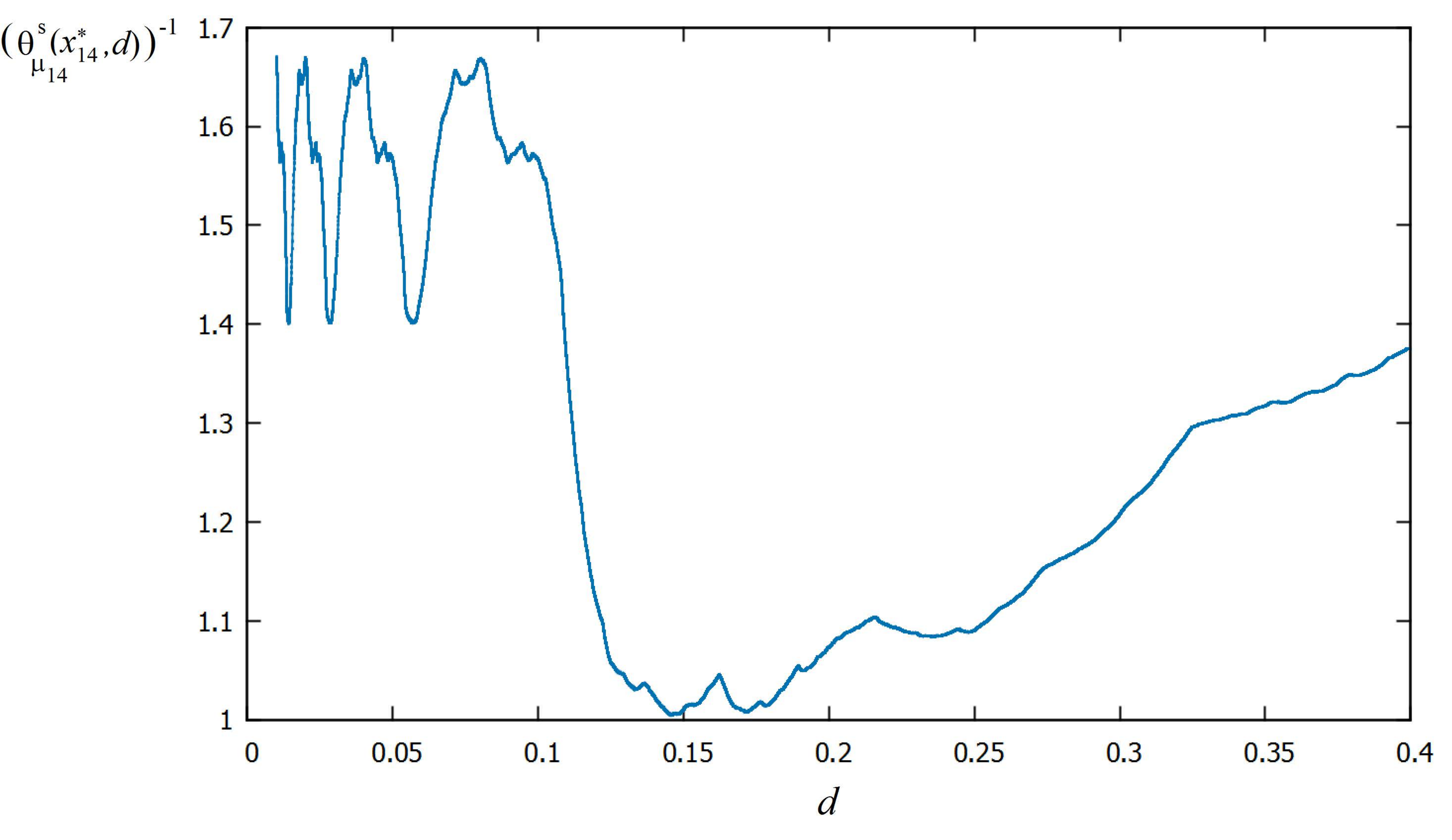} 
\caption{Inverse of the density at the optimal point.\newline
Graph of $(\protect\theta _{\protect\mu _{14}}^{s}(x_{14}^{\ast },d))^{-1}$
as a function of $d$ with $x_{14}^{\ast }=f_{010}(z_{2}).$ The minimal value
of $(\protect\theta _{\protect\mu _{14}}^{s}(x_{14}^{\ast },d))^{-1}$
corresponds to $C^{s}(S)\sim 1.0049$ and the maximal one to $P^{s}(S)\sim
1.6683.$}
\label{invdenbolmaxmin}
\end{figure}

The evolution of the inverse density $(\theta _{\mu
_{14}}^{s}(f_{010}(z_{2}),d))^{-1}$ as a function of the radius, $d,$ is
plotted in Fig.~\ref{invdenbolmaxmin}. Notice that the minimum and maximum
values of this function correspond to the approximations $C^{s}(S)\sim
1.0049 $ and $P^{s}(S)\sim 1.6683,$ respectively.

\section{An upper bound and a conjecture on the spherical Hausdorff measure
of $S$\label{section spherical}}

In Sec.~\ref{subsection metric measures} we have pointed out the analogies
in the definitions of the Hausdorff measure, the spherical Hausdorff measure
and the Hausdorff centred measure. The fundamental difference lies in the
class of sets used to build the coverings. Regarding the spherical Hausdorff
measure, the constraint used in \eqref{premeasure} that the balls be centred
in the set disappears. This easing of restrictions complicates the search
for the optimal covering by adding degrees of freedom to a constraint
optimisation problem already complicated by the natural mismatch between the
rich geometry of a self-similar set and the preeminent regularity of the
Euclidean balls. Nevertheless, a suitable way to deal with the problem of
computability is to use the density function characterisation of this
measure provided by Mor\'{a}n in \cite{Mo1}, which is valid for any
self-similar set, $E,$ satisfying the OSC 
\begin{equation}
\mathcal{H}_{sph}^{s}(E):=\inf \Big\{\frac{(2r)^{s}}{\mu (B(x,r))}:x\in 
\mathbb{R}^{n}\Big\}.  \label{denspherical}
\end{equation}

In the case of the Hausdorff measure, another degree of freedom is added as
there is no constraint on the shape of the covering sets (see Sec.~\ref%
{subsection metric measures}).

Since the first upper bound for the value of $\mathcal{H}^{s}(S),$ obtained
by Marion \cite{MARI} in 1987, and over two decades of work by different
authors (cf. \cite{J}, \cite{JZZ}, \cite{MARI}, \cite{Mora} and \cite{ZF}),
the best bounds for $\mathcal{H}^{s}(S)$ available to date that we
acknowledge are 
\begin{equation}
0.77\leq\mathcal{H}^{s}(S)\leq0.819161232881177,  \label{Moraupbnd}
\end{equation}
given by P. Mora in \cite{Mora}. This improvement on the existing bounds is
obtained through an algorithm based on geometric methods developed by B. Jia
in \cite{J}, and aimed at finding, for each $k\in\mathbb{N}^{+},$ a set of
minimal inverse density among all possible unions of $k$-cylinders.

This method is intrinsically related to the characterisation of the
Hausdorff measure for a self-similar set, $E,$ satisfying the OSC given by
Mor\'{a}n in \cite{Mo1}, 
\begin{equation*}
\mathcal{H}^{s}(E):=\inf\Big\{\frac{|A|^{s}}{\mu(A)}:A\ \text{ is a convex
polytope}\Big \}.
\end{equation*}

Using this same principle, namely, to find a set (which, in this case, can
be an open, convex or closed set) that minimises the density function
mentioned above, Mora constructed in \cite{Mora} a polygon centred at the
barycenter of $T$ (consisting of a collection of $30$-th generation
triangles) whose density provides the upper bound given in \eqref{Moraupbnd}.

In accordance with the geometry of the set obtained in \cite{Mora}, it is to
be expected that the ball optimising \eqref{denspherical}, when $E=S,$ is
also centred at the barycenter of $T.$ Following this conjecture, we have
designed an algorithm to explore the $\mu_{k}$-densities, for $k=14$, of all
the balls centred at the barycentre of $T,$ $x_{b}=(0.5,\frac{\sqrt{3}}{6}), 
$ and with radius $d\in\lbrack\frac{\sqrt{3}}{12},\frac{\sqrt{3}}{3}].$ This
is the range of radii between the radius of the inscribed circle in the
triangle $T\backslash\cup_{i\in M}T_{i}$ and the radius of the circumscribed
circle of $T.$ The aim is to find a ball that minimises the inverse density $%
(\theta_{\mu_{14}}^{s}(x_{b},d))^{-1}.$

In this case, and because the $\mu $-measure of the balls with radii close
to $\frac{\sqrt{3}}{12}$ is almost zero, it is computationally more
convenient to write the algorithm so that it approximates the maximum value
of the density function, 
\begin{equation}
C_{sph}(x_{b}):=\max \Big\{\theta _{\mu }^{s}(x_{b},d):\ \frac{\sqrt{3}}{12}%
\leq d\leq \frac{\sqrt{3}}{3}\Big\},  \label{contdensbari}
\end{equation}%
rather than the minimum value of the inverse density function. The
comparison between the value of $C_{sph}(x_{b})$ and its discrete version%
\begin{align}
C_{sph}^{k}(x_{b})& :=\max \Big\{\theta _{\mu _{k}}^{s}(x_{b},d):\text{ }%
d=|x-x_{b}|,\ x\in A_{k}\Big\}  \label{disdensbar} \\
& =\max \Big\{\theta _{\mu _{k}}^{s}(x_{b},d):\text{ }d=|x-x_{b}|,\
x=(x_{1},x_{2})\in A_{k}\cap S_{2}\text{ and }x_{1}\leq 0.5\Big\}  \notag
\end{align}%
is contained in the following lemma. Note that the last equality holds by
the symmetry of $S,$ and that if $d=|x-x_{b}|$ with $x\in A_{k},$ then $%
\frac{\sqrt{3}}{12}\leq d\leq \frac{\sqrt{3}}{3}.$

\begin{lemma}
\label{bariballdens} Let $C_{sph}(x_{b})$ and $C_{sph}^{k}(x_{b})$ given by %
\eqref{contdensbari} and \eqref{disdensbar}, respectively. \newline
Then, for every $k\geq 3$%
\begin{equation*}
\underline{C}_{sph}^{k}(x_{b})\leq C_{sph}(x_{b})\leq \overline{C}%
_{sph}^{k}(x_{b})
\end{equation*}%
where 
\begin{equation*}
\underline{C}_{sph}^{k}(x_{b}):=\frac{\mu _{k}(B(x_{b},D_{k}-2^{-k}))}{%
(2D_{k})^{s}},\text{ }\overline{C}%
_{sph}^{k}(x_{b}):=K_{k}^{sph}C_{sph}^{k}(x_{b}),
\end{equation*}%
$K_{k}^{sph}:=(1+2^{2-k}\sqrt{3})^{s},$ and $C_{sph}^{k}(x_{b}):=\theta
_{\mu _{k}}^{s}(x_{b},D_{k}).$

\begin{proof}
Let $D\in\big[\frac{\sqrt{3}}{12},\frac{\sqrt{3}}{3}\big]$ be such that $%
C_{sph}(x_{b})=\theta_{\mu}^{s}(x_{b},D).$

Since $D\geq\frac{\sqrt{3}}{3}$ and $k\geq3,$ $f_{01}(z_{2})\in
B(x_{b},D)\cap A_{k}$ and then Lemma \ref{lemma cotas} (ii) can be applied.
Then we get $y_{k}\in A_{k}$ satisfying 
\begin{equation*}
\mu(B(x_{b},D))\leq\mu_{k}(B(x_{b},d_{k}^{\ast})),
\end{equation*}
where $d_{k}^{\ast}=\left\vert y_{k}-x_{b}\right\vert $ and $D-2^{-k}\leq
d_{k}^{\ast}\leq D+2^{-k}.$ Notice that by the symmetry of $S$ we can choose 
$y_{k}$ satisfying the restrictions in (\ref{disdensbar}). Then, 
\begin{align*}
C_{sph}(x_{b}) & =\theta_{\mu}^{s}(x_{b},D)\leq\frac{\mu_{k}(B(x_{b},d_{k}^{%
\ast}))}{(2D)^{s}}=\Big(\frac{d_{k}^{\ast}}{D}\Big)^{s}\theta_{\mu
_{k}}^{s}(x_{b},d_{k}^{\ast})\leq\Big(\frac{D+2^{-k}}{D}\Big)^{s}\theta
_{\mu_{k}}^{s}(x_{b},d_{k}^{\ast}) \\
& \leq\Big(1+\frac{2^{-k}}{D}\Big)^{s}\theta_{\mu_{k}}^{s}(x_{b},D_{k})%
\leq(1+2^{2-k}\sqrt{3})^{s}C_{sph}^{k}(x_{b})=\overline{C}_{sph}^{k}(x_{b}),
\end{align*}
where the last inequality holds true because $D\geq\frac{\sqrt{3}}{12}.$%
\newline
The reverse inequality is immediate using \eqref{contdensbari} and part (i)
of Lemma \ref{lemma cotas}, 
\begin{equation*}
{C}_{sph}(x_{b})=\theta_{\mu}^{s}(x_{b},D)\geq\theta_{\mu}^{s}(x_{b},D_{k})=%
\frac{\mu(B(x_{b},D_{k}))}{(2D_{k})^{s}}\geq\frac{%
\mu_{k}(B(x_{b},D_{k}-2^{-k}))}{(2D_{k})^{s}}=\underline{C}_{sph}^{k}(x_{b}).
\end{equation*}
\end{proof}
\end{lemma}

\begin{figure}[H]
\centering
\includegraphics[scale=0.4]{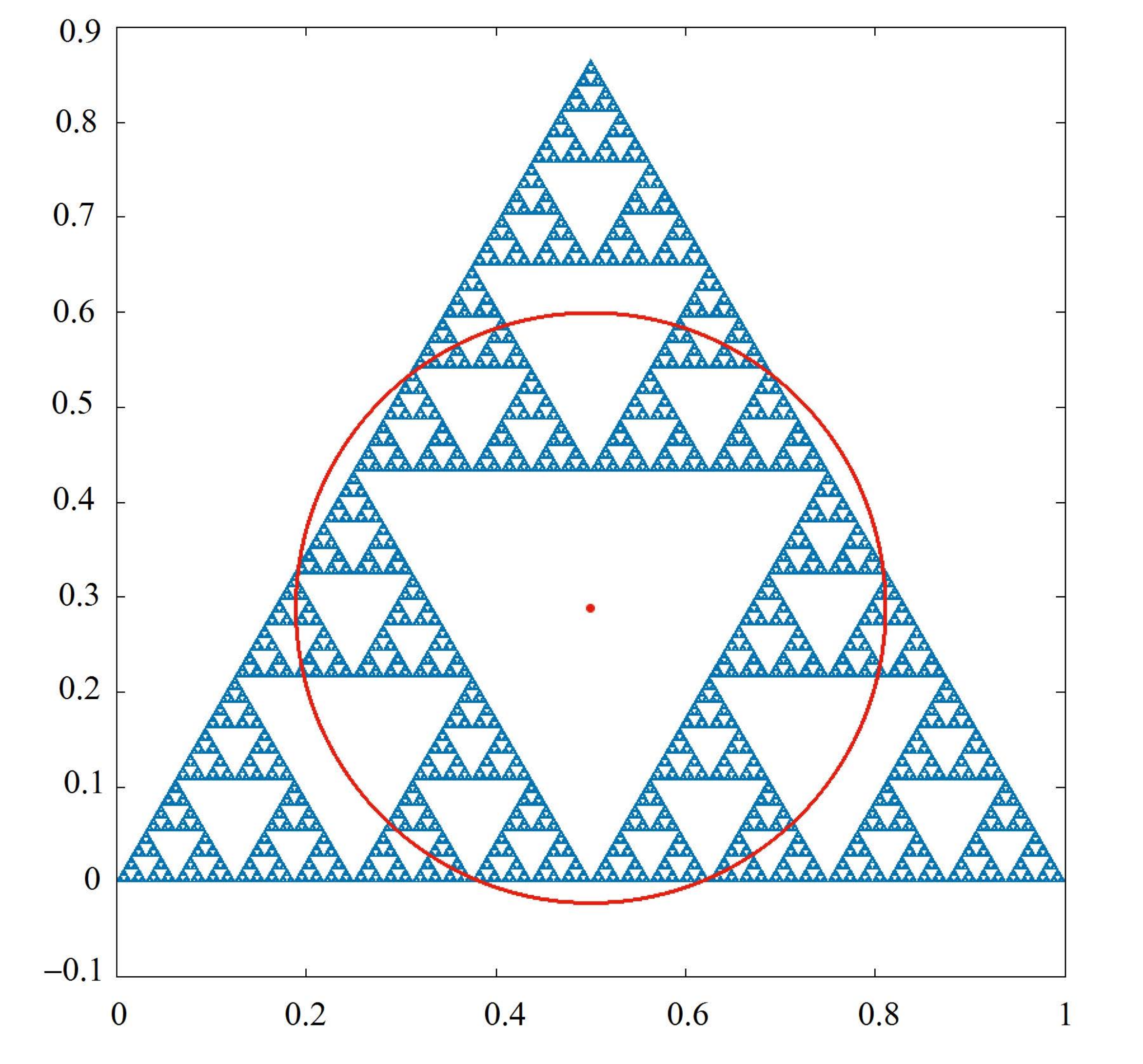} 
\caption{Optimal ball at the barycentre of $T.$\newline
Ball selected by our algorithm when searching for the ball that maximises
the value of $\protect\theta _{\protect\mu _{14}}^{s}(x_{b},d)$ among those
balls centred at the barycentre of $T,$ $x_{b}=(0.5,\frac{\protect\sqrt{3}}{6}),$ and with radius within $[\frac{\protect\sqrt{3}}{12},\frac{\protect
\sqrt{3}}{3}].$ }
\label{bariball}
\end{figure}

\begin{figure}[H]
\centering
\includegraphics[scale=0.4]{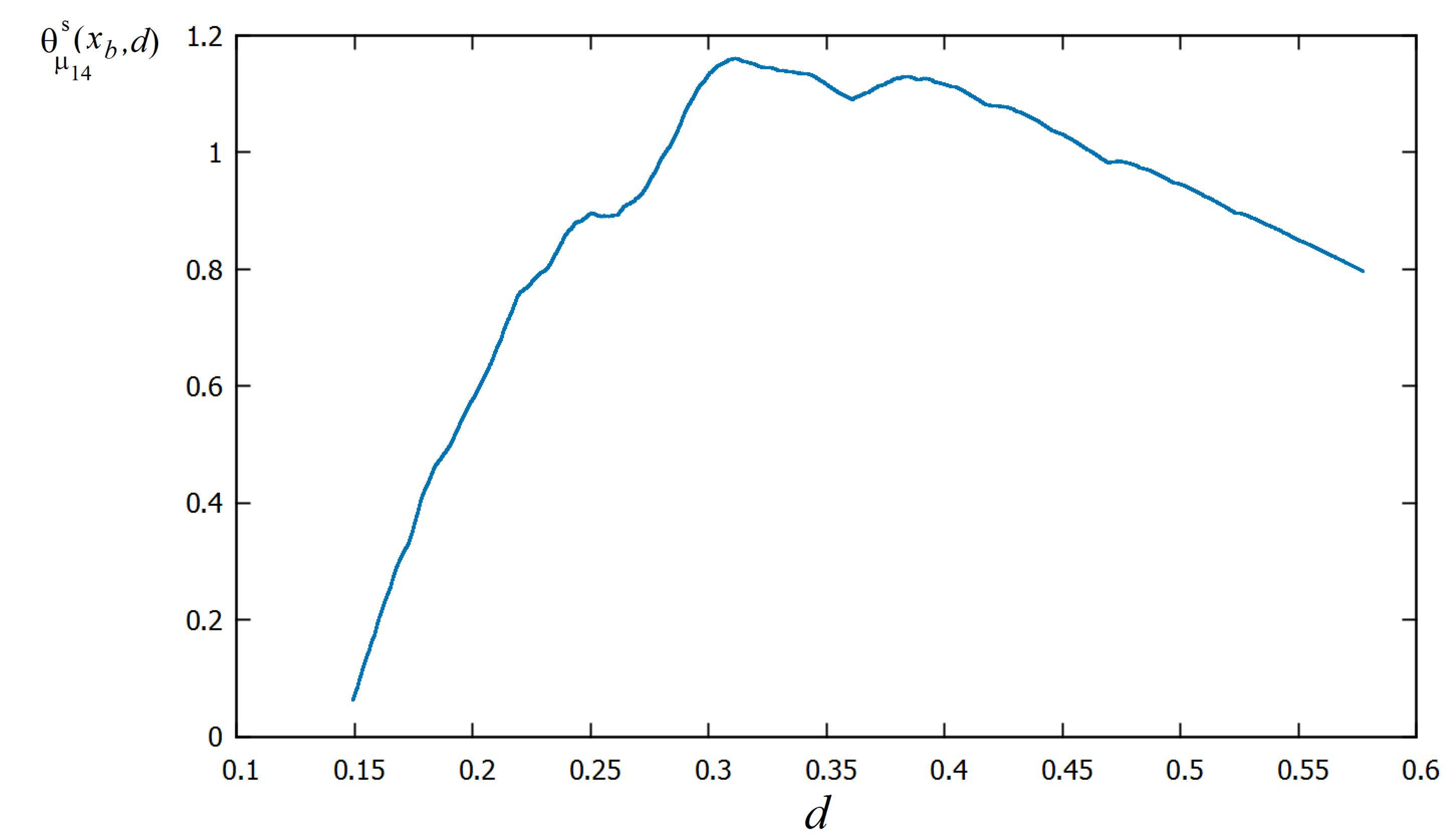} 
\caption{Density function $\protect\theta _{\protect\mu _{14}}^{s}(x_{b},d)$
for $d\in \left[ \frac{\protect\sqrt{3}}{12},\frac{\protect\sqrt{3}}{3}%
\right] .$}
\label{denbari}
\end{figure}

The estimate, $C_{sph}^{k}(x_{b}),$ and bounds, $\underline{C}%
_{sph}^{k}(x_{b})$ and $\overline{C}_{sph}^{k}(x_{b}),$ for $C_{sph}(x_{b})$
obtained by our algorithm for $k=14$ are 
\begin{equation*}
C_{sph}^{14}(x_{b})=1.160630,\text{ }\underline{C}%
_{sph}^{14}(x_{b})=1.160235,\text{ and }\overline{C}%
_{sph}^{14}(x_{b})=1.161408.
\end{equation*}%
These values are rounded to six decimal places, to the smallest value for $%
\underline{C}_{sph}^{14}(x_{b}),$ to the nearest for $C_{sph}^{14}(x_{b}),$
and to the largest for $\overline{C}_{sph}^{14}(x_{b}).$ The radius of the
optimal ball (plotted in Fig.~\ref{bariball}) is $D_{14}\sim 0.3108.$
The maximum value of $\theta _{\mu _{14}}^{s}(x_{b},d),$ $C_{sph}^{14}(x_{b}),$ can be checked in Fig.~\ref{denbari}, where the graph
of the density function $\theta _{\mu _{14}}^{s}(x_{b},d)$ for $d\in \lbrack 
\frac{\sqrt{3}}{12},\frac{\sqrt{3}}{3}]$ is depicted.

The above results, together with Lemma~\ref{bariballdens}, provide the
estimate $\left( C_{sph}^{14}(x_{b})\right) ^{-1}=0.8616$ of $\left(
C_{sph}(x_{b})\right) ^{-1}$ and the bounds%
\begin{equation*}
0.8610\leq (\overline{C}_{sph}^{14}(x_{b}))^{-1}\leq \left(
C_{sph}(x_{b})\right) ^{-1}\leq (\underline{C}_{sph}^{14}(x_{b}))^{-1}\leq
0.8619
\end{equation*}%
which, in turn, allows us to get an upper bound for the spherical Hausdorff
measure of the Sierpinski gasket, namely, 
\begin{equation}
\mathcal{H}_{sph}^{s}(S)\leq \left( C_{sph}(x_{b})\right) ^{-1}\leq 0.8619
\label{sphupbound}
\end{equation}%
and to conjecture that%
\begin{equation*}
\mathcal{H}_{sph}^{s}(S)\sim (C_{sph}^{14}(x_{b}))^{-1}\sim 0.8616.
\end{equation*}%
Observe that our upper bound for $\mathcal{H}_{sph}^{s}(S)$ in (\ref{sphupbound}) is fairly close to the upper bound (\ref{Moraupbnd}) for $%
\mathcal{H}^{s}\mathcal{(S)},$ indicating that the optimal ball found by our
algorithm is a rather efficient covering set for the Sierpinski gasket.

\begin{remark}
Recall (see Remark \ref{invariant measure}) that having an estimate of $%
\alpha (S),$ allows one to estimate $\alpha (A)$ for any Borel set $A,$ and $%
\alpha \in \mathcal{M}^{s}\mathcal{\lfloor }S.$ In particular, taking $%
\alpha =C^{s}\lfloor _{S},$ \ and $A=B(x,d)$ in (\ref{muvshauss}) we get 
\begin{equation}
C^{s}\lfloor _{S}(B(x,d))=C^{s}(S)\mu (B(x,d)).  \label{CsofA}
\end{equation}%
If $k\geq 4,$ and $d>2^{-k},$ then (\ref{CsofA}), Theorem~\ref{theorem
bounds discrete} and Lemma \ref{lemma cotas} (i) guarantee 
\begin{equation}
C_{k}^{\inf }\mu _{k}(B(x,d_{1}))\leq C^{s}\lfloor _{S}(B(x,d))\leq
C_{k}^{\sup }\mu _{k}(B(x,d_{2}))  \label{centball}
\end{equation}%
where $d_{1}:=$ $d-2^{-k}>0,$ and $d_{2}:=d+2^{-k}.$

\begin{example}
\label{cs-of-ball}If we consider $k=14;$ the ball $B(x_{b},D_{14})$ that
maximises (\ref{disdensbar}) as $B(x,d)$ in (\ref{centball}); the estimates 
\begin{equation*}
\mu _{14}(B(x_{b},d_{1}))=0.546105,\text{ }\mu
_{14}(B(x_{b},D_{14}))=0.546290,\text{ \ and \ }\mu
_{14}(B(x_{b},d_{2}))=0.546447,
\end{equation*}%
obtained using our algorithm; together with the estimates of $C_{14}^{\inf },
$ $C_{14}$ and $C_{14}^{\sup }$ of Table~\ref{numresults}, then we get as
estimate of $C^{s}\lfloor _{S}(B(x_{b},D_{14}))$ the value $C_{14}\mu
_{14}(B(x_{b},D_{14}))=0.548968$ and the bounds 
\begin{equation*}
0.547803\leq C^{s}\lfloor _{S}(B(x_{b},D_{14}))\leq 0.549513.
\end{equation*}%
These values are rounded to six decimal places: to the smallest value for
the lower bound, to the largest for the upper bound, and to the nearest for $%
C_{14}\mu _{14}(B(x_{b},D_{14})).$
\end{example}
\end{remark}

\section*{Acknowledgements}

This work was supported by the Universidad Complutense de Madrid and Banco
de Santander (PR108/20-14).

\end{document}